\documentclass[a4paper,11pt]{amsart}

\usepackage{amsfonts,amssymb,amsthm,amsmath}
\usepackage[english]{babel}
\usepackage[all]{xy}
\usepackage{graphics}

1

\newtheorem{teo}{Theorem}[section]
\newtheorem{pro}[teo]{Proposition}
\newtheorem{lemma}[teo]{Lemma}
\newtheorem{cor}[teo]{Corollary}

\theoremstyle{remark}
\newtheorem{oss}{Remark}

\theoremstyle{definition}
\newtheorem{defin}{Definition}

\newcommand{\Z}{\mathbb{Z}}
\newcommand{\N}{\mathbb{N}}

\newcommand{\R}{\mathbb{R}}
\newcommand{\Q}{\mathbb{Q}}
\newcommand{\HHH}{\mathbb{H}}
\newcommand{\G}{\mathbb{G}}
\newcommand{\proj}{\mathbb{P}}

\newcommand{\m}{\mathfrak{m}}
\newcommand{\M}{\mathfrak{M}}

\newcommand{\CCC}{\mathcal{C}}
\newcommand{\AAA}{\mathcal{A}}
\newcommand{\KKK}{\mathcal{K}}
\newcommand{\HH}{\mathcal{H}}
\newcommand{\FF}{\mathcal{F}}
\newcommand{\SSS}{\mathcal{S}}
\newcommand{\OO}{\mathcal{O}}

\DeclareMathOperator{\Tor}{Tor}
\DeclareMathOperator{\lk}{lk}
\DeclareMathOperator{\st}{St}
\DeclareMathOperator{\spec}{Spec}
\DeclareMathOperator{\Hom}{Hom}
\DeclareMathOperator{\projdim}{projdim}
\DeclareMathOperator{\Ass}{Ass}
\DeclareMathOperator{\Ht}{ht}
\DeclareMathOperator{\ann}{ann}
\DeclareMathOperator{\codim}{codim}
\DeclareMathOperator{\cork}{cork}
\DeclareMathOperator{\pos}{pos}

\newcommand{\frec}{\rightarrow}
\newcommand{\iniett}{\hookrightarrow}

\newcommand{\ergo}{\Rightarrow}
\newcommand{\sse}{\Leftrightarrow}

\newcommand{\HR}{\widetilde{H}}  

\newcommand{\fr}{\rightarrow}

\begin{document}

\author{Silvano Baggio}
\title[equivariant K-theory of smooth toric varieties]{equivariant K-theory \\of smooth toric varieties }
\address{Dipartimento di Matematica\\
Universit\`a di Bologna\\
Piazza di Porta San Donato 5\\
40126 Bologna\\ Italy}
\email{baggio@dm.unibo.it}

\begin{abstract}We characterize the smooth toric varieties for which the Merkurjev spectral sequence, connecting equivariant and ordinary K-theory, degenerates. We find under which conditions on the support of the fan the $E^2$ terms of the spectral sequence are invariants by subdivisions of the fan. Assuming these conditions, we describe explicitly the $E^2$ terms, linking them to the reduced homology of the fan.

\mbox{}

\end{abstract}
\maketitle

\tableofcontents


\section*{Introduction}
  The aim of this work is to highlight some relationships between the equivariant K-theory of a smooth toric variety and the combinatorics of the associated fan. 

We start from a result due to Merkurjev (\cite{mer}), that gives a comparison between the equivariant and the ordinary K-theory: if $X$ is a smooth toric variety, there is a spectral sequence: 
$$
E_{pq}^2=\textrm{Tor}^{RT}_p(K_q^T(X),\Z)\Longrightarrow K_{p+q}(X).
$$

 Our first problem  (Section 2) is to characterize the class of smooth toric varieties for which Merkurjev's spectral sequence degenerates at $E^2$, that is, such that $E^2_{pq}=0$ for all $p\not=0$.  It was known that this class includes all complete toric varieties, which correspond to fans whose support covers the whole space $\R^n$. Their K-theory can be immediately expressed in terms of the equivariant K-theory: $K_q(X)=K_q^T(X)\otimes_{RT}\Z$. 

 Vezzosi and Vistoli (\cite[\S 6]{vevi}) give two descriptions of the equivariant K-theory of a smooth toric variety $X(\Delta)$, 

(1) as a subring of the product $\prod_{\sigma\in\Delta_{\max}}K_*(k)\otimes RT_{\sigma}$, and 

(2) by means of a presentation similar to the Stanley-Reisner algebra over the simplicial complex $S_{\Delta}$, associated to the fan $\Delta$ (The vertices of $S_{\Delta}$ correspond to the 1-dimensional cones of the fan, and the simplexes correpond to the cones).
 Both descriptions have been extremely useful. 

 Proposition \ref{piattezza} states that Merkurjev's spectral sequence degenerates if and only if the following conditions hold:
\begin{itemize}
\item[H1.] $\HR_i(\lk\sigma,\Z)=0$ for all $\sigma\in\Delta$ and $i<\dim(\lk\sigma)$, where $\HR$ is the reduced simplicial homology, and $\lk\sigma$ is the {\em link} of $\sigma$, that is, the set of the faces, disjoint from $\sigma$, of the faces in $S_{\Delta}$ that contain $\sigma$; 
\item[H2.] $\HR_i(S_{\Delta},\Z)=0$ for all $i<\dim(S_{\Delta})$.
\end{itemize}

 This is an easy consequence of a well-known theorem by Reisner: the vanishing of the reduced homology (in all degrees but the top one) of a simplicial complex and of all of its links is equivalent to saying that its Stanley-Reisner ring is Cohen-Macaulay. However, assuming the (independent) Theorem \ref{principale}, we can prove Proposition \ref{piattezza} without reference to Reisner's theorem.

 Conditions H1 and H2 involve only {\em the support} $|\Delta|$ of
 $\Delta$: they are independent of the way it is subdivided into
 cones. Indeed, they depend only on the topology of $|\Delta|$, in
 particular  H1 is equivalent to saying that $S_{\Delta}$ is a
 homology manifold. A natural question is, whether it is possible to find a description of the terms $E_{pq}^2$ also depending only on the topology of $|\Delta|$. The answer (Section 3) is negative for all toric varieties whose fan does not satisfy the local condition H1: this is a straightforward consequence of Proposition \ref{ostruzione}.  

Therefore, in Section 4 we {\em assume} that condition H1 holds, and moreover that all the maximal cones of $\Delta$ have the same dimension as the variety. In other words, (a geometric realization of) $S_{\Delta}$ is a submanifold (with boundary) of the sphere $S^{n-1}\subset\R^n$.
 Under these assumptions, we can find an explicit expression of the terms $E_{pq}^2$, where only the reduced homology of $S_{\Delta}$ appears. In fact we have (Theorem \ref{principale}):
 
\begin{equation*}
\Tor^{RT}_p(K_0^T(X(\Delta)),\Z)  \cong  \bigoplus_{i=p+1}^n\HR^{^{i-p-1}}\Bigl(S_{\Delta},\Z^{^{\binom{n}{i}}}\Bigr).
\end{equation*}
The terms $\Tor_p^{RT}(K_q^T(X),\Z)$, for $q>0$, can be obtained from those with $q=0$, applying the Theorem of Universal Coefficients, since (by the Stanley-Reisner presentation) $K_q^T(X)=K_0^T(X)\otimes_{\Z}K_q(k)$.

 The proof of Theorem \ref{principale} is based on the idea (already used in \cite{bri},\cite{BBFK} and \cite{brelu}) to give $\Delta$ the topology induced by the face order, and consider suitable sheaves of rings on it. In our case, the expression of  $K_0^T(X)$ as a subalgebra of $\prod_iRT_{\sigma_i}$ suggests to construct a sheaf $\AAA$, whose ring of global sections is $K_0^T(X)$. The terms $E^2_{p,0}$ in Merkurjev's spectral sequence are the hypercohomology groups of the Koszul complex of $\AAA$. This hypercohomology is itself the limit of another spectral sequence, but, thanks to condition H1, most of its $E_1$ terms vanish, and we can draw useful information from it.        

\subsection*{Acknowledgments} I would like to thank my thesis advisor, Angelo Vi\-stoli, who suggested me the fascinating subject of toric varieties, and helped me considerably during the making of my PhD thesis (whose results are exposed in this paper). I also thank Michel Brion for his careful reading it and for giving me interesting comments.   


\section{Preliminaries on toric varieties and K-theory}
\subsection{Toric varieties}
The most complete references for toric varieties are Fulton's and
Oda's books, \cite{F} and \cite{oda}. More information can be found in  Danilov's \cite{dan} and Brylinski's \cite{bry} papers. For a more detailed introduction, see also \cite{gvt}.
 We recall here some terminology and properties.

We fix a field $k$. Let $\G_m$ be the multiplicative group of $k$. For a positive integer $n$, $T=T^n=(\G_m)^n$ is the n-dimensional (algebraic) \emph{torus} on $k$. 
  
 A {\em toric variety} is a normal algebraic variety $X$ on $k$, with a dense embedding $T\iniett X$, such that the action of $T$ on itself by multiplication can be extended to an action of $T$ on $X$.

Let $N\cong \Z^n$ be a lattice, embedded in $N_{\R}=\R\otimes N\cong
\R^n$. A {\em (rational covex polyhedral) cone} $\sigma$ is the convex
hull of a finite set of vectors $v_i\in N$:
 \begin{equation*}
 \sigma = \pos\{v_1,\dots,v_s\} \stackrel{\text{def}}{=}
 \{r_1v_1+\dots+r_sv_s\mid r_i\geq 0\}.
 \end{equation*}

 Let $M=\Hom(N,Z)$. The {\em dual cone} of the cone $\sigma$ is $\check{\sigma}=\{u\in M_{\R}\mid
 \langle u, v\rangle\ge 0\,\forall v\in \sigma\}$. A {\em support hyperplane} for $\sigma$ is the set
$\{v\in N\mid \langle u, v\rangle=0\}$,  for some $u\in \check{\sigma}$.  A {\em face} $\tau$ of $\sigma$ is the intersection of $\sigma$ and a support hyperplane; we write $\tau\prec\sigma$. 

If $\sigma$ is a cone in $N_{\R}$, the semigroup $S_{\sigma}=\check{\sigma}\cap M$ is finitely generated (\cite[Proposition 1]{F}). Let $A_{\sigma}=k[S_{\sigma}]$ the semigroup algebra of $A_{\sigma}$. The {\em affine toric variety associated to $\sigma$} is (the affine algebraic variety) $X_{\sigma}\stackrel{\text{def}}{=}\spec(A_{\sigma})$.

\begin{defin} A {\em fan} of cones in $N_{\R}\cong\R^n$ is a finite set $\Delta$ of cones in $N_{\R}$, with the following properties:
\begin{enumerate}
 \item if $\sigma\in\Delta$ and $\tau\prec\sigma$, then $\tau\in \Delta$,
 \item if $\sigma_1,\sigma_2\in\Delta$, $\sigma_1\cap\sigma_2$ is a common face of $\sigma_1$ and $\sigma_2$.
 \end{enumerate}
The {\em support} of a fan $\Delta$ is the union of its cones, $|\Delta|=\bigcup\{\sigma\in\Delta\}\subset N_{\R}$. We set $\Delta_i=\{\sigma\in \Delta\mid\dim\sigma=i\}$, while   $\Delta_{\max}$ is the set of the {\em maximal} cones, that is, cones that are faces of no other cones in $\Delta$.
 \end{defin}
Given a fan $\Delta$, there is a toric variety $X(\Delta)$ associated
to it: it is obtained glueing together the affine toric varieties
corresponding to its cones. See  \cite[Ch.~1]{F} for
details. The lattices $N$ and $M$ are, respectively, the groups of 1-parameter
subgroups and the group of characters of $T$ (\cite[1.2]{oda}). Conversely, every toric variety can be obtained from a fan,
via this construction (see \cite{bry}). 

$X(\Delta)$ is \emph{smooth} if and only if every $\sigma\in\Delta$ is generated by a subset of some basis of the lattice $N$ \cite[2.1]{F}: we say then that $\Delta$ is {\em regular}. In particular a regular fan is {\em simplicial}, that is, every cone admits a number of generators equal to its dimension.
 
Let $N$ and $N'$ be two lattices, $\Delta$ a fan in $N_{\R}$ and
$\Delta'$ in $N'_{\R}$, $\psi:N\to N'$ a homomorphism, and
$\varphi=\psi\otimes \R:N_{\R}\to N'_{\R}$. If, for each
$\sigma\in\Delta$, $\varphi(\sigma)\subset\sigma'$ for some
$\sigma'\in\Delta'$, then $\varphi$ induces a {\em toric morphism},
that is, a morphism of algebraic varieties $X(\Delta)\to X(\Delta')$, equivariant with respect to the actions of the dense tori. Such a toric morphism is proper if and only if  $\varphi^{-1}(|\Delta'|)=|\Delta|$ (see \cite[2.4]{F}).  $X(\Delta)$ is {\em complete} if and only if  $|\Delta|=N_{\R}\cong\R^n$.

 \subsection{Orbits and their closures}(\cite[3.1]{F})\label{orbite}
There is a 1-1 correspondence between orbits for the action of $T$ and cones in $\Delta$: for each cone $\sigma$ of dimension $d$ the corresponding orbit $\mathcal{O}_{\sigma}$ has dimension $n-d$. Moreover:
 $$
 \tau\prec\sigma\sse\mathcal{O}_{\sigma}\subset\overline{\mathcal{O}}_{\tau}
 $$
The closure $\overline{\OO}_{\sigma}$ of the orbit $\OO_{\sigma}\subset X(\Delta)$ is the finite union of smaller orbits, corresponding to the cones in the {\em star} $\st\sigma=\{\tau\in\Delta\mid\sigma\prec\tau\}$. Moreover, $\overline{\OO}_{\sigma}$ is itself a toric variety, for the action of the torus  $T/T_{\sigma}$. In fact, if $N_{\sigma}=N\cap\langle\sigma\rangle$ is the lattice of one parameter subgroups of $T_{\sigma}$, let us consider the quotient lattice $N(\sigma)=N/N_{\sigma}$, and the real vector space $N(\sigma)_{\R}=N_{\R}/(N_{\sigma})_{\R}$. The projection $\bar{\tau}$ on $N(\sigma)_{\R}$ of a cone $\tau$ in $N_{\R}$ is also a cone, and the set $\{\bar{\tau}\subset N(\sigma)_{\R}\mid \tau\in\st\sigma\}$ is a fan in $N(\sigma)_{\R}$. The associated toric variety is isomorphic to the orbit closure $\overline{\OO}_{\sigma}$.
\subsection{The $RT$-algebras $RT_{\sigma}$}\label{RT}
The stabilizer of the points in $\mathcal{O}_{\sigma}$ is a subtorus
$T_{\sigma}\subset T$, whose lattice of one parameter subgroups is
$N\cap\langle\sigma\rangle\subseteq N$. Let
$M=\widehat{T}=\textrm{Hom}(T,\G_{\textrm{m}})$ be the group of
characters of $T$, and  $RT=\Z \widehat{T}=\Z[M]$ the representation
algebra of $T$. For each cone $\sigma$,
$RT_{\sigma}=\Z[M/\sigma^{\perp}]$, where $\sigma^{\perp}=\{m\in M\mid
\langle u, m\rangle=0\,\forall\: u\in\sigma\}$ is the orthogonal subgroup
in $M$ to $\sigma$ (or to the subspace generated by $\sigma$). If $\tau\prec\sigma$, the inclusion induces a homomorphism $RT_{\sigma}\rightarrow RT_{\tau}$. In particular, for each $\sigma$, the inclusion $T_{\sigma}\subset T$ induces a morphism $RT\to RT_{\sigma}$, so each $RT_\sigma$ is a $RT$-algebra.  

We can easily get the projective dimension of $RT_{\sigma}$ as a $RT$-module, indeed we can give explicit free resolutions of each  $RT_{\sigma}$. 
 
If $\sigma$ is a cone of dimension $d$, $\sigma^{\perp}$ is
a subgroup of rank $n-d$ in $M$. Let $m_1,\dots,m_{n-d}$ be  a basis of $\sigma^{\perp}$.   
 Since the elements $m_1-1,\dots,m_{n-d}-1\in RT$ generate (as an ideal) the kernel of the projection $RT\to RT_{\sigma}$ (that is, $\Z[M]\to\Z[M/\sigma^{\perp}]$), and they are  a regular sequence in $RT$, the associated Koszul complex (\cite[\S 16]{mat}),
 \begin{equation}\label{Koszul}
 K.(\sigma^{\perp})=K.(m_1-1,\dots,m_{n-d}-1)
 \end{equation}
is acyclic, and $H_0(K.)=RT_{\sigma}$: it is a free resolution of $RT_{\sigma}$. It follows that   $\projdim_{RT}(RT_{\sigma})=n-\dim(\sigma)$. Moreover, if we consider $\Z$ as a $RT$-module via the rank map $RT\to \Z$, and we tensor $K.$ with  $-\otimes_{RT}\Z$, the differentials vanish; 
 so we have
 \begin{equation}\label{TorRT}
 \Tor^{RT}_i(RT_{\sigma},\Z)=H_i(\Z\otimes K.)\cong\begin{cases}\bigwedge^i \sigma^{\perp}\cong \Z^{^{\binom{n-d}{i}}} & \text{for $0\leq i\leq n-d$}\\0 & \text{for $i>n-d$.}\end{cases}
 \end{equation}
%
 
 \subsection{Two results in the K-theory of toric varieties}
An introduction to the K-functor in algebraic geometry can be found in (\cite{man}). Quillen  (\cite{qui}) defined higher K-theory groups, as homotopy groups of the geometric realization of a category. 

If $X$ is a scheme (or an algebraic space) with the action of an algebraic group, equivariant coherent sheaves can be defined, and the corresponding K-theory is called equivariant K-theory: the basic definitions and theorems are in \cite{tho}.  

  In equivariant K-theory computations are often easier than in ordinary K-theory, and information about ordinary K-theory can be recovered from equivariant K-theory, for example by means of spectral sequences.   

In fact a recent theorem by Vezzosi and Vistoli expresses the equivariant K-theory of a smooth toric variety as a subalgebra of a product of representation rings of tori, while the comparison of the equivariant and the ordinary K-theory is provided by a spectral sequence, introduced by A.~Merkurjev (from results of M.~Levine). 

From now on, given a smooth toric variety $X$ over the field $k$, with the action of an algebraic torus $T$, we will denote by $K_*(X)$ the ordinary K-theory, and with $K_*^T(X)$ the $T$-equivariant K-theory. 

Notice that $K_0^T(X)$ is a module over the representation ring $RT\cong K_0^T(\spec k)$, via the homomorphism induced by the morphism $X\to \spec k$ (\cite[Example 2.1]{mer}).

 \begin{teo}[\textmd{\cite[Theorem 6.2]{vevi}}] 
\label{immersione}
If $X=X(\Delta)$ is a smooth toric variety, and  $\sigma_1,\dots,\sigma_r$ are the maximal cones of  $\Delta$, there is an injective homomorphism of $RT$-algebras
 $$
 K_*^T(X)\hookrightarrow \bigoplus_{i=1}^r K_*(k)\otimes
 RT_{\sigma_i}.
 $$

An element $(a_i)\in\bigoplus_{i=1}^r K_*(k)\otimes RT_{\sigma_i}$
is in the image of this homomorphism if and only if, for each  $i\neq  j$ the restrictions of  $a_i$ and $a_j$ to $K_*(k)\otimes RT_{\sigma_i\cap\sigma_j}$ coincide.
 \end{teo}

 \begin{teo}[\textmd{\cite[Theorem 4.3]{mer}}]
Let $X$ be a smooth toric variety, with maximal torus $T$. There is a homology spectral sequence
 $$
 E_{pq}^2=\textrm{Tor}^{RT}_p(K_q^T(X),\Z)\Longrightarrow K_{p+q}(X),
 $$
such that the boundary homomorphisms
 $$
 \Z\otimes_{RT}K_q^T(X)\longrightarrow K_q(X)
 $$
are induced by the functor that forgets the action of $T$ $K_*^T(X)\rightarrow K_*(X)$.

In particular the ring homomorphism $\Z\otimes_{RT}K_0^T(X)\rightarrow K_0(X)$ is an isomorphism.
\end{teo}


\section{Conditions on the fan $\Delta$ for the $RT$-flatness of $K_0^T(X(\Delta))$}
In \cite{vevi} Vezzosi and Vistoli give, as a consequence of their description of the equivariant K-theory, a sufficient condition for the Merkurjev spectral sequence to degenerate at $E^2$: in that case there is an isomorphism $\Z\otimes_{RT}K^T_*(X)\cong K_*(X)$. The condition is that the action of the big torus {\em admits enough limits} (\cite[Definition 5.8]{vevi}). A smooth toric variety $X(\Delta)$ admits enough limits if and only if the set
\begin{equation}\label{abbast}
\bigcap_{\tau\in\Delta}\bigcup_{\sigma\in\st\tau}(\sigma+\langle\tau\rangle)\subset N_{\R}
\end{equation}
has nonempty interior (\cite[Proposition 6.7]{vevi}).  This fact is very useful, as it allows to single out many toric varieties whose K-theory has an explicit description. For instance, it follows immediately from the above condition that complete toric varieties belong to this class.

 Nonetheless, there are two main flaws in it. First, it provides only a sufficient, but not a necessary condition: we will see below an example of a toric variety that does not have enough limits, but such that the Merkurjev sequence degenerates. Second, it involves not only the face order of the cones in the fan, but the exact way in which the cones are set inside the vector space $N_{\R}$.  

In this section we will give a necessary and sufficient condition, depending only on the simplicial structure of the fan. This condition follows rather easily from the presentation of $K^T_*(X)$ as a Stanley-Reisner ring (or ``face ring'') of the simplicial complex $S_{\Delta}$ associated to $\Delta$, and from a well-known criterion by Reisner, giving conditions on a simplicial complex in order for its Reisner-Stanley ring to be Cohen-Macaulay. However, in section \ref{alt}, after we prove Theorem \ref{principale}, we will be able to avoid using Reisner's Theorem: Proposition \ref{piattezza} will be a corollary to Theorem \ref{principale} and  Proposition \ref{present}.\vspace{5pt}

 We recall here some notions of combinatorics, which we use in Proposiotion \ref{piattezza}.

By a {\em simplicial complex} over the finite set (of {\em vertices}) $V=\{v_1\dots,v_s\}$ we mean the pair $(V,\Sigma)$, where $\Sigma$ is a set of subsets of $V$ (the {\em simplexes} or {\em faces}), such that:
\begin{align*}
 &\sigma\in\Sigma,\, \tau\subset\sigma\ergo\tau\in\Sigma; \\
 &\{v\}\in\Sigma\quad\forall\, v\in V.
\end{align*}
The {\em dimension} of a face $\sigma$ is the number of vertices of $\sigma$ minus one; the dimension of  $\Sigma$ is  $\max\{\dim\sigma\mid\sigma\in\Sigma\}$. A simplicial complex is {\em pure} if all its maximal faces have the same dimension.

Let $\Delta$ be a regular fan in $\R^n$. We associate to it the {\em simplicial complex} $S_{\Delta}$ in the following  way:
\begin{itemize}
\item[-] the {\em vertices} of $S_{\Delta}$ correspond to one-dimensional cones of $\Delta$; 
\item[-] the {\em faces} correspond to cones in $\Delta$.
\end{itemize}
The dimension of $S_{\Delta}$ is $\max\{\dim\sigma\mid\sigma\in\Delta\}-1$. A geometric realization of $S_{\Delta}$ is $|\Delta|\cap S^{n-1}\subset\R^n$, where $S^{n-1}=\{x\in\R^n\mid|x|=1\}$ .

\begin{defin}
 Let $R$ be a ring, and $\Sigma$ a simplicial complex with vertices $V=\{v_1,\dots,v_s\}$ and faces $\SSS\subset\mathcal{P}(V)$. The {\em Stanley-Reisner algebra} on $R$ relative to $\Sigma$ is the $R$-algebra
$$
R[\Sigma]=\frac{R[X_1,\dots,X_s]}{I},
$$ 
where $I\subset R[X_1,\dots,X_s]$ is the ideal generated by all monomials $X_{i_1}\cdots X_{i_h}$ with $\{v_{i_1},\dots,v_{i_h}\}\not\in \SSS$.
\end{defin}
\begin{defin}(See \cite[5.3]{bh})
Let $(V,\Sigma)$ be a simplicial complex of dimension $n-1$, and let $V$ be given a total order. For each $i$-dimensional face $\sigma$ we write $\sigma=[v_0,\dots,v_i]$ if $\sigma=\{v_0,\dots,v_i\}$ and $v_0<\dots<v_i$.

The {\em augmented chain complex} of $\Sigma$ is:  
$$
\CCC(\Sigma):0\to\CCC_{n-1}\stackrel{\partial}{\to}\CCC_{n-2}\to\dots\to\CCC_0\stackrel{\partial}{\to}\CCC_{-1}\to 0
$$ 
where we set 
$$
\CCC_i=\bigoplus_{\stackrel{\sigma\in\Sigma}{\dim\sigma=i}}\Z\sigma\quad\text{and}\quad \partial\sigma=\sum_{j=0}^i(-1)^j\sigma_j
$$
for all $\sigma\in\Sigma$, and $\sigma_j=[v_0,\dots,\widehat{v_j},\dots,v_i]$ for $\sigma=[v_0,\dots,v_i]$. By definition, $\dim\emptyset=-1$.

For an abelian group $G$, the $i$-th {\em reduced simplicial homology} of $\Sigma$ with values in $G$ is:
$$
\HR_i(\Sigma,G)=H_i(\CCC(\Sigma)\otimes G)\quad i=-1,\dots,n-1.
$$
The dual (cochain) complex $\Hom_{\Z}(\CCC(\Sigma),G)$ has differentials $\bar{\partial}$, defined as: $(\bar{\partial}\phi)(\alpha)=\phi(\partial \alpha)$, for  $\phi\in\Hom_{\Z}(\CCC_i,G)$, $\alpha\in\CCC_{i+1}$.   
The $i$-th group of {\em reduced simplicial cohomology} of $\Sigma$ with values in $G$ is:
$$
\HR_i(\Sigma,G)=H_i(\Hom_{\Z}(\CCC(\Sigma),G))\quad i=-1,\dots,n-1.
$$
\end{defin}
If $\sigma$ is a face of the simplicial complex $\Sigma$, the {\em link} of $\sigma$ in $\Sigma$ is $\lk_{\Sigma}\sigma=\lk\sigma=\{\tau\in\Sigma\mid \sigma\not\subset\tau,\,\sigma\cup\tau\in\Sigma\}$. It is easy to see that $\lk_{\Sigma}\sigma$ is itself a simplicial complex over the set $\{v\in V\mid v\in\tau\;\exists\,\tau\in\lk\sigma\}$.
\begin{pro}[Reisner's Criterion]\label{reisn}(\cite[Theorem 3]{rei})
Let $(V,\Sigma)$ be a simplicial complex. The ring  $\Z[\Sigma]$ is  Cohen-Macaulay if and only if, for every cone  $\sigma\in\Sigma$ we have
$$
\HR_i(\lk\sigma,\Z)=0 \quad \forall i<\dim (\lk\sigma),
$$
and, moreover,
$$
\HR_i(S_{\Delta},\Z)=0 \quad  \forall i<\dim (S_{\Delta}).
$$
\end{pro}

The following proposition gives a presentation of the equivariant K-theory of smooth toric varieties, in terms of generators and relations, similar to that given in \cite{BDP} for the equivariant cohomology.
\begin{pro}(\cite[Proposition 6.4]{vevi})\label{present}

There is an isomorphism of $K_*(k)$-algebras
$$
K_*^T(X_{\Delta})\cong\frac{K_*(k)[X_{\rho}^{\pm 1}]_{\rho\in\Delta_1}}{I},
$$
where $X_{\rho}$ are indeterminates, each corresponding to a one-dimensional cone of $\Delta$, and $I$ is the ideal in $K_*(k)[X_{\rho}^{\pm 1}]$ generated by all products $\prod_{\rho\in F}(X_{\rho}-1)$, with $F$ a subset of $\Delta_1$ that does {\em not} generates a cone of $\Delta$.
\end{pro}
Now we can state and prove the result we announced.
\begin{pro}\label{piattezza}
Let $X=X(\Delta)$ a smooth $n$-dimensional toric variety, associated to the fan $\Delta$. The following are equivalent:
\begin{enumerate}
\item $K_0^T(X)$ is a flat ($\Leftrightarrow$ projective) $RT$-module;
\item $Tor_i^{RT}(K_0^T(X),\Z)=0$ for all $i>0$;
\item $Tor_i^{RT}(K_q^T(X),\Z)=0$ for all $q\geq 0$ and for all $i>0$;
\item The following two conditions hold:
\begin{itemize}
\item[(i)] $\Delta_{\max}=\Delta_n$  ($\ergo S_{\Delta}$ is pure); 
\item[(ii)]
$
\HR_i(\lk\sigma,\Z)=0 \quad \forall i<\dim (\lk\sigma),\;\forall\sigma\in S_{\Delta},
$
 and 
$$
\HR_i(S_{\Delta},\Z)=0 \quad  \forall i<\dim (S_{\Delta}).
$$
\end{itemize}
\end{enumerate}
\end{pro}
\begin{proof}
Let us write for short  $K=K_0^T(X)$, and set $\Delta_{\max}=\{\sigma_1,\dots,\sigma_r\}$. We will show $(4)\sse(1)\sse(2)\sse(3)$.\vspace{5pt}\\
$(4)\Rightarrow (1)$ If $\dim(\sigma)=n$, $RT_{\sigma}=RT$, so condition  $(i)$ and  Theorem  \ref{immersione} imply that $K$ embeds into the product $RT^r$. The ring $RT$ is a subring of  $K$ via the diagonal map ($RT\iniett RT^r$, $a\mapsto(a,\dots,a)$), and $K$ is a finitely generated $RT$-module.\\
Note also that, if $\M$ is a maximal ideal in $K$, and $\m=\M\cap RT$, then $\dim(K_{\M})=\dim(RT_{\m})$. In fact, if $\wp$ is a minimal prime ideal in $K$, contained in $\M$, then $K/\wp\iniett RT/\bar{\wp}$ for some minimal prime ideal $\bar{\wp}$ in $RT$, such that $\bar{\wp}\subset\m$.\\ 
We have verified the hypotheses of  the following lemma.
\begin{lemma}\label{ohoh}
Let $R$ and $S$ two noetherian rings, $R$ regular, $R\subset S$, and $S$ a finite $R$-module. Moreover, suppose that for all ideals  $\M$ maximal in $S$ and $\m$ maximal in $R$, with $\M\cap R=\m$, we have $\dim R_{\m}=\dim S_{\M}$.

 Then $S$ is a projective $R$-module if and only if it is a Cohen-Macaulay ring.
\end{lemma}
\begin{proof}  \emph{(Lemma \ref{ohoh})} Let $\m$ a maximal ideal in  $R$. By the Auslander-Buchs\-baum formula (\cite[Thm 1.3.3]{bh}), $S_{\m}$ is a projective $R_{m}$-module if and only if the depth of $S_{\m}$ as a $R_{\m}$-module equals the dimension (and so the depth, as $R$ is regular) of the ring $R_{\m}$. 

 If $d=\dim R_{\m}$ and $\bar{x}=x_1,\dots,x_d$ is a regular system of parameters for  $R_{\m}$, then it is a regular sequence  in $S_{\m}$, and this is true if and only if $\bar{x}$ is a regular sequence in  $S_{\M}$ for each maximal ideal $\M$ in $S$ such that $\M\cap R=\m$. The sequence $\bar{x}$ is also maximal, for each choice of such an $\M$, by the assumption that $\dim S_{\M}=\dim R_{\m}$.
\end{proof}
 
It remains to prove that $K$ is a Cohen-Macaulay ring.  

 Proposition \ref{present} gives a ring isomorphism  $K\cong \Z[X^{\pm 1}]/I$, with $X=X_1,\dots,X_l$ indeterminates, and $I=(\prod_{j\in F}(X_j-1)|F\not\in S_{\Delta})\subset \Z[X^{\pm 1}]$.

On the other hand, the Stanley-Reisner ring of $S_{\Delta}$ is $\Z[Y]/J$, with $Y=Y_1,\dots,Y_l$, and $J=(\prod_{i_j\in F}Y_{i_j}|F\not\in S_{\Delta})\subset \Z[Y]$.

Consider the following ring homomorphism: 
\begin{eqnarray}
\label{phi}\Z[Y]&\stackrel{\varphi}{\to}&\Z[X]\\\notag
 Y_i & \mapsto & X_i-1
\end{eqnarray}
Obviously  $\varphi$ is an isomorphism, and  $\varphi(J)=I$, so $\varphi$ induces an isomorphism
$$
\frac{\Z[Y]}{J}\cong\frac{\Z[X]}{I},
$$
and it remains an isomorphism if we localize respectively at the multiplicative systems  $S_J=\{\prod_i Y_i^k\}_{k\in\Z}$ and $S_I=\{\prod_i(X_i-1)^k\}_{k\in\Z}$:
$$
\frac{S_J^{-1}\Z[Y]}{J}\cong\frac{\Z[X^{\pm 1}]}{I}=\frac{S_I^{-1}\Z[X]}{I}.
$$
The ring on the left is a localization of the Stanley-Reisner ring $\Z[S_{\Delta}]$, so it is Cohen-Macaulay if the latter is such: this is true by Proposition \ref{reisn}.\vspace{5pt}

$(1)\Rightarrow (4)$ From what we said above it follows that $(1)\ergo(4\,.ii)$, that is: $\Z[S_{\Delta}]$ is Cohen-Macaulay. We need only to show that the maximal cones have all dimension $n=\dim(X)$.

Let us prove this by contradiction: suppose that, say, $\dim\sigma_1<n$.

This implies that, if we denote by $\pi_i:RT\rightarrow RT_{\sigma_i}$ the projections induced by $T_{\sigma_i}\subset T$, and set $I_i=\ker \pi_1$, we have: $\Ht I_1\ge 1$. 

Let  $\pi^{i}_{ij}$ be the projections $RT_{\sigma_i}\to RT_{\sigma_i\cap\sigma_j}$. By Theorem \ref{immersione},\\ $K=\{(x_1,\dots,x_r)\in\prod_{i=1}^r RT_{\sigma_i}\mid\pi_{ij}^i(x_i)=\pi_{ij}^j(x_j) \quad \forall i,j\}$.

If $x_1$ is a nonzero element in $\bigcap_{i>1}\ker\pi_{1i}^1$, then $x=(x_1,0,\dots,0)\in K$, and $\ann x=I_1$, so $\Ht(\ann x)\ge 1$.

Let $\m$ be the kernel of the rank map $RT\rightarrow \Z$,  $A$ the localization  $RT_{\m}$, and $\wp$ an associated prime ideal of $K$ as a $RT$-module, such that $I_1\subset\wp$. As $\wp\subset\m$, $\wp A\in \Ass_A K_{\m}$ too, and $\Ht(\wp A)=\Ht(\wp)$.

Since, for a local noetherian ring  $R$ and a finitely generated $R$-module $M$, $\projdim_R M\geq \Ht(\mathfrak{a})$ for any ideal $\mathfrak{a}$ associated of $M$, we have: 
$$1\le \Ht(\wp A)\le d\stackrel{\text{def}}{=}\projdim_A(K_{\m}).$$
 It follows that $0\not=\Tor_d^A(K_{\m},\Q)=A\otimes_{RT}\Tor_d^{RT}(K,\Z)$, so $\Tor_d^{RT}(K,\Z)\not=0$, contradicting our assumption.\vspace{5pt}

$(1)\Rightarrow (2)$ Obvious.\vspace{5pt}

$(2)\Rightarrow (1)$ 
From now on, set $A=\Z[Y_1,\dots,Y_n]$. Notice that the morphism $\varphi$ defined in \eqref{phi}, followed by the canonical arrow to the localization at the multiplicative system $S_I=\{\prod_i(X_i-1)^k\}_{k\in\Z}$, is a ring monomorphism that  makes $RT$ a flat $A$-module. Remember that  $\forall\,\sigma\in\Delta$, $RT_{\sigma}\cong RT/I_{\sigma}$ for some ideal $I_{\sigma}$ contained in the kernel of $RT\to\Z$. Let $J_{\sigma}=I_{\sigma}\cap A$, and $A_{\sigma}=A/I_{\sigma}$. For any $\tau\prec\sigma$ we have a surjective map $A_{\sigma}\to A_{\tau}$ such that $RT_{\sigma}\to RT_{\tau}$ is obtained from it, by tensoring with $-\otimes_A RT$. We define here 
$$
M=\Big\{(a_1,\dots,a_r)\in\prod_{i=1}^r A_{\sigma_i}\:\Big|\: {a_i}_{|A_{\sigma_i\cap\sigma_j}}={a_j}_{|A_{\sigma_i\cap\sigma_j}}\:\forall\,i,j\Big\}.
$$    
Then $K=M\otimes_{A} RT$, and $K$ is flat over  $RT$ if and only if  $M$ is flat, or projective, over $A$. Assumption (2) implies that 
\begin{equation}\label{tora}
\Tor_i^A(M,\Z)=0\quad\forall\,i>0.
\end{equation}
As $A$ is a graded ring with $A_0=\Z$, $M$ is flat over $A$ if and only if $M_p$ is flat over $A_p$, with $p\in\N$ varying among all primes, where we have defined  $I_p=(p)+(Y_1,\dots,Y_n)$, $A_p=A_{I_p}$ and $M_p=M_{I_p}=M\otimes_A A_p$.

By the local criterion of flatness (see for ex. \cite[Theorem 6.8]{eis}), $M_p$ is flat over $A_p$ if and only if $\Tor^{A_p}_1(M_p,\Z/(p))=0$. From \eqref{tora} it follows that  $\Tor_1^{A_p}(M_p,\Z_{(p)})=0$, as $\Z_{I_p}\cong\Z_{(p)}$, the localization of $\Z$ at $p$. From the short exact sequence 
$$
0\to\Z_{(p)}\stackrel{\cdot p}{\to}\Z_{(p)}\to\Z/(p)\to 0
$$  
it follows that $\Tor_1^{A_p}(M_p,\Z/(p))=0$ if and only if multiplication by $p$ in $\Z_p\otimes_{A_p}M_p$ is injective: and this is true since  $M$, and also $M_p$, are torsion free. \vspace{5pt}

$(3) \ergo (2)$ Obvious. \vspace{5pt}

$(2) \ergo (3)$ Fix an integer $q>0$. Proposition \ref{present} implies that 
\begin{displaymath}
K_q^T(X)=K_q(k)\otimes_{\Z} K_0^T(X).
\end{displaymath}  
 Let $L_*$ be the Koszul resolution \eqref{Koszul} of $\Z$: every $L_i\cong\bigwedge^i (RT)^n$ is a free $RT$-module. Then, by the definition of  $\Tor^{RT}_*(-,\Z)$ as a derived functor of  $-\otimes_{RT}\Z$,    
\begin{eqnarray*}
\Tor_p^{RT}(K_0^T(X),\Z) & = & H_p(L_*\otimes_{RT} K_0^T(X));\\
\Tor_p^{RT}(K_q^T(X),\Z) & = & H_p(L_*\otimes_{RT} K_q^T(X))=\\
                         & = &H_p(L_*\otimes_{RT} (K_0^T(X)\otimes_{\Z}K_q(k))).
\end{eqnarray*}
Since $L_p\otimes_{RT} (K_0^T(X)\otimes_{\Z}K_q(k))\cong (L_p\otimes_{RT} K_0^T(X))\otimes_{\Z}K_q(k)$ we can apply the Universal Coefficient Theorem (in homology: see for example \cite[V.11, Theorem 11.1]{mac}): given a homology complex $L_*$ of abelian groups with no elements of finite order, and an abelian group $G$, then for any $i>0$ there is a splitting short exact sequence:
$$
0\frec H_i(L_*)\otimes G\frec H_i(L_*\otimes_{\Z} G)\frec \Tor^{\Z}_1(H_{i-1}(L_*),G)\frec 0.
$$
In our case, for every  $p>0$, there is a short exact sequence  ($K=K_0^T(X)$):
\begin{multline}\label{kappasup}
0\frec \Tor^{RT}_p(K,\Z)\otimes_{\Z} K_q(k)\frec \Tor^{RT}_p(K^T_q(X),\Z)\frec \\ \frec \Tor_1^{\Z}(\Tor_{p-1}^{RT}(K,\Z),K_q(k))\frec 0.
\end{multline}
To conclude, notice that, for  $p>1$, $\Tor^{RT}_p(K_0(X),\Z)=0$, while $\Z\otimes_{RT} K_0(X)=\Tor_0^{RT}(K_0(X),\Z)$ is torsion-free, so  $\Tor_1^{\Z}(\Z\otimes K_0(X),G)=0$ for any abelian group $G$.
\end{proof}
\textbf{Example} of a toric variety that satisfies the conditions of  Proposition \ref{piattezza}, but that does not admit enough limits. \vspace{5pt}

Let $\{e_1,e_2,e_3\}$ a basis in $N\cong\Z^3$. For any triple $(a_1,a_2,a_3)\in\{0,1\}^3$, let  $\sigma_{a_1,a_2,a_3}$ be the cone generated by  $\{(-1)^{a_1}e_1,(-1)^{a_2}e_2,(-1)^{a_3}e_3\}$ (they are  the ``octants'' of $\R^3$).

Let us consider now the fan $\Delta$, such that the maximal cones are:
$$
\Delta_{\max}=\{\sigma_{000},\sigma_{010},\sigma_{011},\sigma_{111}\}.
$$
$X(\Delta)$ is a non-complete variety; it can be embedded as an open $T$-invariant set of  $\proj^1\times\proj^1\times\proj^1$ (more precisely: it can be obtained from  $\proj^1\times\proj^1\times\proj^1$ removing the closure of three 1-dimensional orbits).

Let $\rho_1=\pos\{e_2\}$ and $\rho_2=\pos\{-e_1\}$.
 
Notice that $\sigma_{000}$ is the only maximal cone that contains $\rho_1$ as a face, while $\sigma_{111}$ is the only maximal cone containing  $\rho_2$. So 
$$\bigcup\limits_{\sigma\in\st\rho_1}\sigma+\langle\rho_1\rangle=\sigma_{000}\cup\sigma_{010},$$ 
while 
$$\bigcup\limits_{\sigma\in\st\rho_2}\sigma+\langle\rho_2\rangle=\sigma_{011}\cup\sigma_{111}.$$
The intersection of these two sets is  $\pos\{e_1,-e_2\}$, which has empty interior. A fortiori the set  \eqref{abbast} has empty interior, so $X(\Delta)$ does not admit enough limits. 

On the other hand, $|\Delta|\cap S^2$ is homeomorphic to   $\{(x,y)\in\R^2\mid x^2+y^2\leq 1\}$, so it is contractible, therefore the reduced homology of  $S_{\Delta}$ and of every link is zero: the equivalent conditions of Proposition  \ref{piattezza} are thus satisfied.


\section{K-theory and subdivisions of the fan} 
The combinatorial conditions equivalent to the flatness of $K_0^T(X(\Delta))$, given in Proposition \ref{piattezza}, depend only on the topology of the support of the fan $\Delta$, not on its subdivision into cones. If $K_0^T(X)$ is not flat, some $E^2$ term outside the column of index zero in Merkurjev spectral sequence is nonzero: 
$$\Tor^{RT}_p(K^T_q(X),\Z)\not=0\quad\exists\, q,\;\exists\, p>0.$$ 
A natural question is: given the toric variety $X(\Delta)$, does there exist a description of these groups $\Tor^{RT}_p(K^T_q(X(\Delta)),\Z)$, involving only topological invariants of $|\Delta|$? We will show in Proposition  \ref{ostruzione} and in Corollary \ref{ostruzione2} that a necessary and sufficient condition for the existence of such a description is the vanishing of the reduced homology of all links of the cones $\sigma\not=0$ in $\Delta$, in every dimension strictly less than the dimension of $\lk \sigma$.

 In the next section we will assume this condition on $\Delta$, and will be able to give the expected description.  

\subsection{Blow-ups along orbit closures}
Let us begin by comparing two toric varieties associated to fans with the same support. Let $\Delta_1$ and $\Delta_2$ be two fans in $\R^n$, corresponding to toric varieties  $X_1=X(\Delta_1)$ and $X_2=X(\Delta_2)$. Assume that $|\Delta_1|=|\Delta_2|$.
 Let $\Delta'$ be a common subdivision of the two fans: that is, every cone of $\Delta_1$ or $\Delta_2$ is union of cones of $\Delta'$. 
The identity of $\R^n$ induces two $T$-equivariant morphisms $\phi_i:X(\Delta')\frec X_i$.
 They are {\em birational}, as the two varieties have the same dimension, so they contains the same torus $T^n$ as a dense open subset, and  {\em proper} (see \cite[2.4]{F}).  

For birational and proper equivariant morphisms between toric varieties the weak decomposition theorem holds; it was proved independently by  R.~Morelli in \cite{mor} and by J.~W\l odarczyk in \cite{wlo}.
\begin{teo}\label{decomposizione}Every birational proper equivariant morphism between smooth toric varieties $\phi:X(\Delta')\to X(\Delta)$ can be decomposed as follows:
$$
X(\Delta')\to X_1\gets X_2\to \dots \gets X_k \to X(\Delta), 
$$
where $X_1,\dots,X_k$ are smooth toric varieties, and every arrow is an equivariant morphism  obtained by composing blow-ups along $T$-invariant closed subvarieties. 
\end{teo}
 We have reduced the problem to the study of  $K^T_0$ and its groups $Tor^{RT}_*$, when $\phi$ is the blow-up of an orbit closure.

 Blowing up the orbit closure relative to the cone  $\sigma$ consists of replacing $\bar{\OO}_{\sigma}$ with a $T$-invariant divisor. This corresponds to the modification of the fan, called {\em star subdivision}. 

Given two cones  $\sigma,\tau$ in $N_{\R}$, let us denote with $\sigma+\tau$ their {\em Minkowski sum}: $\sigma+\tau=\{x+y\mid x\in\sigma,\,y\in\tau\}$. It is a cone in $N_{\R}$. 
\begin{defin} Let $\Delta$ be an arbitrary fan in $N_{\R}$, $\sigma\in\Delta$, and $\rho\subset\sigma$ a 1-dimensional cone (ray) passing through the relative interior of $\sigma$ (i.e. $\rho$ is contained in no proper face of $\sigma$). 
The  {\em star subdivision} of $\Delta$, relative to the cone $\sigma$ and to the ray $\rho$, is the fan $\{\sigma/\rho\}\cdot\Delta$, obtained as the union of $\Delta\setminus\st\sigma$ and
$$
\{\rho+\sigma'+\nu\mid\sigma\not=\sigma'\prec\sigma,\,\nu\in\lk\sigma\}.
$$
\end{defin}
In other words, this subdivision does not change the cones not containing $\sigma$, but it divides each cone containing $\sigma$ into the cones generated by $\rho$ and every proper face.

Assume now that $\sigma$ is regular. Then $\sigma$ is generated by vectors in $N$, $e_1,\dots,e_s$, that can be completed to a basis of $N$. let $\rho_{\sigma}$ the 1-dimensional cone generated by $q=e_1+\dots+e_s$. Let us denote with $\{\sigma\}\cdot\Delta=\{\sigma/\rho_{\sigma}\}\cdot\Delta$.
\begin{pro}(\cite[Prop. 1.25]{oda})
The blow-up of the toric variety $X(\Delta)$ along the orbit closure $\bar{\OO}_{\sigma}$ is the toric variety associated to the fan obtained by the star subdivision of  $\Delta$ relative to $\sigma$ and $\rho_{\sigma}$:
$$
Bl_{\bar{\OO}}X(\Delta)=X_{\{\sigma\}\cdot\Delta}
$$ 
\end{pro}

\subsection{How the groups $\Tor^{RT}_p(K_0^T(X(\Delta)),\Z)$ change by subdivisions of $\Delta$}
Let us see how the equivariant K-theory of the blow up depends on the K-theory of the variety, of the orbit closure and of the exceptional divisor.  
Let  $Y=\bar{\OO}_{\sigma}$ be the centre of the blow up, and  $Y'=\phi^{-1}(Y)$ the exceptional divisor. The situation can be pictured by the following (cartesian) square:
$$
\xymatrix{Y'\ar@{^{(}->}[r]^j \ar@{-{>>}}[d]_{\psi} & X'\ar@{-{>>}}[d]^{\phi} \\
          Y \ar@{^{(}->}[r]^i        & X
}  
$$
where $i$ and $j$ are the embedding of $Y$ and $Y'$ respectively, and $\psi=\phi_{|Y'}$. $Y'$ is the projective bundle over $Y$ defined by the conormal sheaf $\mathcal{C}$ of $Y$ in $X$: $Y'=\proj(\mathcal{C})$ (See \cite[Theorem II.8.24]{har}). 
 The $K^T$-theory of $X'$ is connected to that of $X$, $Y$ and $Y'$ by a short exact sequence of $RT$-modules:
$$
0\frec K^T_*(Y)\stackrel{\alpha}{\longrightarrow} K_*^T(X)\times K^T_*(Y')\stackrel{\beta}{\longrightarrow} K^T_*(X')\frec 0
$$   
where $\alpha$ and $\beta$ are defined as follows. Let $F$ be the class in $K^T_*(Y')$ of the kernel of the canonical surjection $\psi_*\mathcal{C}\to\mathcal{O}_{Y'}(1)=\mathcal{O}_{\proj(\mathcal{C})}(1)$. If $y\in K^T_*(Y)$, we define $\alpha(y)=(i_*y,-\psi(y)\lambda^{-1}F)$. If $x\in K^T_*(X)$ and $y'\in K^T_*(Y')$, then $\beta(x,y')=\phi^*x-j_*y'$ (The corresponding short exact sequence for ordinary K-theory can be found in \cite[Theorem 15.2]{man}; the equivariant version can be obtained in a similar manner) 

This sequence splits, so the long exact sequence of the $\Tor_*^{RT}(-,\Z)$ itself splits into short exact sequences, for any $i\ge 0$:
\begin{multline}\label{torscop}
0\frec \Tor^{RT}_i(K^T_*(Y),\Z)\frec \Tor^{RT}_i(K_*^T(X),\Z)\times \Tor^{RT}_i(K^T_*(Y'),\Z)\frec\\\frec \Tor^{RT}_i(K^T_*(X'),\Z)\frec 0
\end{multline}

\begin{oss}In the case of toric varieties, we can see that the exceptional divisor of a blow up is a projective bundle, also by comparing the fan $\Delta$ with the star subdivision $\{\sigma\}\cdot\Delta$, or better by comparing the fans associated to $Y$ and $Y'$.

The orbit closure relative to $\sigma$ is the toric variety associated to the projection $\overline{\st\sigma}$ of $\st_{_{\Delta}}\sigma$ on the quotient $N(\sigma)_{\R}$, and that $Y'$ is the orbit closure associated to $\rho_{\sigma}$ in $\{\sigma\}\cdot\Delta$. It is the toric variety associated to the projection $\overline{\st \rho_{\sigma}}$ of $\{\sigma\}\cdot\st\sigma$ on the  quotient $V_{\rho}:=N_{\R}/\langle \rho_{\sigma}\rangle$ (see section \ref{orbite}). 

Let us recall the following definition (\cite[Ch.~VI, Definition 6.3]{ewa}).
\begin{defin}
Let $\Sigma=\Sigma'\cdot\Sigma''=\{\sigma'+\sigma''\mid\sigma'\in\Sigma',\sigma''\in\Sigma''\}$ a fan in $\R^n$ ({\em join} of two fans $\Sigma'$ and $\Sigma''$) such that 
\begin{itemize}
\item[(a)]$\Sigma'$ is contained in a vector subspace $U\subset\R^n$ with $\dim(U)<n$,
\item[(b)]$\Sigma''$ can be projected bijectively on a fan $\Sigma_0$ contained in the orthogonal complement $U^{\perp}$ of $U$.
\end{itemize}
Then we call $\Sigma_0$ a  {\em projection fan} of  $\Sigma$ {\em perpendicular} to $\Sigma'$, and say that $\Sigma$ has a projection fan (with respect to $\Sigma'$, $\Sigma''$).  
\end{defin}
Now denote with $\Delta_{\sigma}$ the projection on $V_{\rho}$ of the set of cones (in $N_{\R}$) $\{\tau+\rho_{\sigma}\mid\sigma\not=\tau\prec\sigma\}$, and with $\Delta_{\lk\sigma}$ the projection, on the same space, of $\{\tau+\rho_{\sigma}\mid \sigma\not\prec\tau\in\Delta,\,\tau+\sigma\in\Delta\}$. We can easily see that $\overline{\st \rho_{\sigma}}=\Delta_{\sigma}\cdot\Delta_{\lk\sigma}$, and that $\Delta_{\lk\sigma}$ can be bijectively projected on the orthogonal complement of $|\Delta_{\sigma}|$. Therefore $\overline{\st \rho_{\sigma}}$ has a projection fan: this is significant in view of the following proposition  (\cite[VI, Theorem 6.7]{ewa}):
\begin{pro}
Let $\Sigma$, $\Sigma'$, $\Sigma''$ be regular fans in $\R^n$, such that $\Sigma=\Sigma'\cdot\Sigma''$, and let $\Sigma_0$ be the projection fan of $\Sigma$ perpendicular to $\Sigma'$. Then the projection $\pi:\Sigma\to\Sigma_0$ induces a map of fans such that, for any $\sigma_0\in\Sigma_0$, we have an isomorphism
$
\bar{\pi}^{-1}(X_{\sigma_0})\cong X(\Sigma')\times X_{\sigma_0}.
$
\end{pro}
This means exactly that $X(\Sigma)\to X(\Sigma_0)$ is a fibre bundle on $X(\Sigma_0)$, with fibre $X(\Sigma')$. In the case of the blow up of $\overline{\OO}_{\sigma}\cong X(\Sigma_0)$, $\Sigma'=\Delta_{\sigma}$, and $X(\Sigma')\cong \proj^{r-1}$.
\end{oss}
 \vspace{5pt} 
The following proposition and corollary contain the condition that we must assume, in order to describe $\Tor^{RT}_p(K_q^T(X),\Z)$ in terms of topological invariants of $|\Delta|$.
\begin{pro}\label{ostruzione}
Let $X=X(\Delta)$ the smooth toric variety associated to the fan $\Delta$ in $N_{\R}\cong \R^n$,  $\sigma\in\Delta_d$  a cone of dimension $d\geq 2$, $\Delta'=\{\sigma\}\cdot\Delta$, $X'=X(\Delta')$  (i.e. $X'$ is the blow up of $X$ along the closure of the orbit $\OO_{\sigma}$). The following hold:
\begin{itemize}
\item[(i)] if $\HR_j(\lk\sigma,\Z)\not=0$ for some $j<d-1$, then $\Tor^{RT}_i(K_q^T(X),\Z)\not\cong\Tor^{RT}_i(K_q^T(X'),\Z)\:\exists i>0\;\forall\,q\geq 0$;
\item[(ii)] conversely, if $\Tor^{RT}_i(K_q^T(X),\Z)\not\cong\Tor^{RT}_i(K_q^T(X'),\Z)\:\exists\,i,q\geq 0$, then for some $\tau\prec\sigma$ we have $\HR_j(\lk\tau,\Z)\not=0$ $\exists j<\dim(\tau)-1$. 
\end{itemize}
\end{pro}
\begin{proof}{\em (i)} With the above remarks and notations in mind,
  we need only to recall the relation between the K-theory of $Y$ and
  $Y'$. We know that $Y'$ is a projective bundle with base $Y$ and
  fibre $\proj^{d-1}$; $T$ acts on the base and on each fibre, and the
  projection $Y'\to Y$ is $T$-equivariant. We apply \cite[Theorem
  3.1]{tho} and get a group isomorphism $K^T_q(Y')\cong  K^T_q(Y)^{\oplus d}$. Hence the abelian group $\Tor^{RT}_i(K_q^T(Y'),\Z)$ is the $d$-fold direct power of $\Tor^{RT}_i(K_q^T(Y),\Z)$. Assume that $\HR_j(\lk\sigma,\Z)\not=0$ for some $j<d-1$. From Proposition \ref{piattezza} it follows that $\Tor^{RT}_i(K_q^T(Y),\Z)\not=0$ for some $i>0$ and for any $q\geq 0$. From the exact sequence \eqref{torscop} we have 
\begin{eqnarray*}
\Tor^{RT}_i(K_q^T(X'),\Z) & \cong & \Tor^{RT}_i(K_q^T(X),\Z)\oplus\Tor^{RT}_i(K_q^T(Y),\Z)^{d-1}\\
&\not\cong& \Tor^{RT}_i(K_q^T(X),\Z).
\end{eqnarray*}
{\em(ii)} By the same exact sequence \eqref{torscop} we can say that, if $\Tor^{RT}_i(K_q^T(X'),\Z)\not\cong\Tor^{RT}_i(K_q^T(X),\Z)$, then $\Tor^{RT}_i(K_q^T(Y),\Z)\not=0$, and, by Proposition \ref{piattezza}, this is true if and only if there exists some $\tau\prec\sigma$ such that $\HR_j(\lk\tau,\Z)\not=0$ $\exists j<\dim(\tau)-1$.   
\end{proof}

\begin{cor}\label{ostruzione2}
Let $\Delta$ be a regular fan in $\R^n$ and $X=X(\Delta)$ the associated smooth toric variety. Then the following are equivalent:
\begin{itemize} 
\item[(a)] for every subdivision  $\Delta'$ of $\Delta$, 
$$\Tor^{RT}_p(K_q^T(X),\Z)\cong\Tor^{RT}_p(K_q^T(X(\Delta')),\Z)$$ 
for any $p>0$ and $q\geq 0$;
\item[(b)] $\HR_j(\lk\sigma,\Z)=0$ for any $\sigma\in S_{\Delta}$ and for any $j<\dim(\lk\sigma)$.
\end{itemize}
\end{cor}  
\begin{proof}
(a)$\ergo$(b) follows immediately from Proposition \ref{ostruzione}.
 
(b)$\ergo$(a) follows from Proposition \ref{ostruzione}, from Theorem \ref{decomposizione}, and from the following remark. 
\end{proof}
\begin{oss}Let $x\in|\Delta|$, and let $\HR_{\bullet}(|\Delta|,|\Delta|\smallsetminus x;\Z)$ be the relative reduced (singular) homology. 
Then condition (b) of Corollary \ref{ostruzione} is equivalent to 
$$
\HR_i(|\Delta|,|\Delta|\smallsetminus x;\Z)\quad\forall i<\dim(|\Delta|),\;\forall x\in|\Delta|,
$$
see e.g. \cite[Proposition 4.3]{sta}: therefore that condition depends only on the topology of $|\Delta|$, and if it holds, it still holds when we perform a star subdivision $\Delta$ or, when possible, the inverse operation (star reunion of cones).
\end{oss}


\section{Computing $\Tor ^{RT}_p(K_q^T(X),\Z)$ -- when possible}

The main result in this section is Theorem \ref{principale}, where we express the groups $\Tor ^{RT}_p(K_0^T(X),\Z)$ in terms of the reduced simplicial homology of the fan $\Delta$.  

The proof of that theorem is based on the following remark: the embedding of $K_0^T(X(\Delta))$ in the product of rings of representations  of tori (Theorem \ref{immersione}), can be restated by saying that $K_0^T(X(\Delta))$ is the ring of global sections of a sheaf of $RT$-algebras over a fan space $\Delta$. 

We will first recall the definition and the main properties of such sheaves. Then we will apply them to compute the groups $\Tor ^{RT}_p(K_0^T(X),\Z)$, under the assumption that $\Delta$ satisfies condition \eqref{nullhom}, i.e. the  vanishing of the local homology of $S_{\Delta}$. By means of homological algebra machinery we will be able to find explicit formulas linking these groups to the reduced homology of $S_{\Delta}$.

\subsection{Cohomology of sheaves on fan spaces}
Some of the following notations and ideas have already been used in \cite{bac}, \cite{brelu}, \cite{BBFK}, \cite{bri} . In \cite{bac} some aspects of sheaves on partially ordered sets are studied, but the main focus is on geometric lattices.  Bressler and Lunts (\cite[\S 3]{brelu})  study sheaves of $\R$-algebras on fan spaces, especially in the non rational case. In \cite{BBFK} Barthel et al. give conditions on the topology of $|\Delta|$, in order for the group of global sections of a sheaf on $\Delta$ to be free on the ring of polynomials. Brion (\cite{bri}) applies the theory of sheaves on fan spaces to the study of the polytope algebra.  

Our reference for definitions and general properties of sheaves on topological spaces is \cite[II.1 and III.1-2]{har}. 
\subsubsection{Sheaves on partially ordered sets}

Let $(X,\leq)$ be a finite partially ordered set. There is a topology on $X$, induced by the order, such that the open sets are the {\em increasing} subsets, that is:
$$
A\subset X,\, A\text{ is open}\sse(x\in A,\,y\in X,\,x\leq y\ergo y\in A).
$$ 
In fact, it is straightforward to verify that the intersection and the union of a family of increasing subsets are increasing subsets. The continuous maps between two partially ordered sets, with this topology, are the order-preserving maps.

For each element $z\in X$, the following sets are defined:
\begin{itemize}
\item[] $\bar{X}_z=\{x\in X: z\leq x\}$ is the smallest open set containing the element $z\in X$;
\item[] $X_z=\{x\in X: z<x\}=\bar{X}_z\smallsetminus \{z\}$;   
\item[] $\bar{X}^z=\{x\in X: x\leq z\}$ is the closure of $\{z\}$;
\item[] $X^z=\{x\in X: x<z\}=\bar{X}^z\smallsetminus \{z\}$.
\end{itemize}
Let us consider now a {\em presheaf} of abelian groups $\FF$ on $X$. To every open set (increasing subset) $U\subset X$ an abelian group $\FF(U)$ is associated: $\FF(U)$ is called the set of {\em sections} of $\FF$ on $U$; for every inclusion of open sets $V\subset U$ there is a group homomorphism $\rho^U_V:U\frec V$, called the {\em restriction} from $U$ to $V$. We will adopt the usual shorthand $s_{|V}=\rho^U_V(s)\quad\forall s\in\FF(U)$.

$\FF$ is a  {\em sheaf} if the following glueing condition for sections holds: 
\begin{quote}
given an open subset $U\subset X$, an open covering $U=\cup_i U_i$, and a family of sections $(s_i)$,  $s_i\in \FF(U_i)\;\forall i$, that are  compatible on the intersections, that is,  ${s_i}_{|U_i\cap U_j}={s_j}_{|U_i\cap U_j}\;\forall i,j$; {\em there exists} a {\em unique} global section $s\in\FF(U)$  such that $s_{|U_i}=s_i\:\forall i$.
\end{quote}    
The {\em stalk} of $\FF$ at $x\in X$ is by definition $\FF_x=\varprojlim\FF(U)$, with $U$ varying in the set $\{U\subset X\mid x\in U,\, U\text{ open}\}$. Since this set contains  $\bar{X}_x$ as a minimum, we have more simply $\FF_x=\FF(\bar{X}_x)$.
\begin{oss}
 We can write in a more explicit form the sections of a sheaf $\FF$ on the open set $U$. Since 
$$U=\bigcup_{x\in U}\bar{X}_x=\bigcup_{x\in U_{\min}}\bar{X}_x,$$
 where $U_{\min}$ is the set of the minimal elements of $U$, we have
\begin{equation}\label{sezioni}
\FF(U)\cong\{(s_x)_{x\in U_{\min}}\mid s_x\in\FF_x,\,{s_x}_{|\bar{X}_x\cap\bar{X}_y}={s_y}_{|\bar{X}_x\cap\bar{X}_y}\;\forall x,y\in U_{\min}\}.
\end{equation}
Notice that, by the above formula, $\FF$ is {\em determined} by the data of \\(1) the stalks $\FF_x$ for every $x\in X$, and \\(2) the restrictions $\FF_x\frec\FF_y$ for all pairs $x\leq y$.

In fact, in order to assign a sheaf on the finite partially ordered set $X$, it is enough to assign an abelian group $G_x$ to every element $x\in X$, and, to each inclusion $x\leq y$, a group homomorphism $\rho_{xy}:G_x\to G_y$, such that for any triple of elements $x\leq y\leq z$ in $X$ the equality $\rho_{yz}\rho_{xy}=\rho_{xz}$ holds. If, given these data, we {\em define} $\FF$ by means of \eqref{sezioni} for every open $U\subset X$, we get a sheaf on $X$, whose stalks are the groups $\{G_x\}$. \vspace{5pt}
\end{oss} 
We recall that a sheaf $\FF$ is {\em flabby} if, for each pair of open set $U\subset V$, the restriction $\FF(V)\frec\FF(U)$ is surjective. Flabby sheaves are acyclic, i.e. they have zero cohomology in positive degree, see \cite{har}, Prop. 2.5: we will exploit this property both to compute the cohomology of the simple sheaves on fan spaces (see below), and to connect the groups  $\Tor^{RT}_*(K^T(X),\Z)$ to the hypercohomology of suitable complexes of sheaves. 

The following criterion for flabbiness is proved in \cite{brelu} in the case $X$ is a fan space.
\begin{lemma}\label{fiacc}
Let $\FF$ a sheaf (of abelian groups) on the partially ordered set $(X,\leq)$. The following are equivalent:
\begin{itemize}
\item[(i)] $\FF$ is flabby,
\item[(ii)] for any $x\in X$, the restriction $\FF_x=\FF(\bar{X}_x)\frec\FF(X_x)$ is surjective. 
\end{itemize}
\end{lemma}
\begin{proof}
Clearly $(i)\ergo (ii)$. Conversely, to prove $(ii)\ergo (i)$, it is
enough to show that, if condition $(ii)$ holds, then, for any open $U$
and any $x\not\in U$, we can extend every section on $U$ to a section
on $V=U\cup \bar{X}_x$. Any subset $U\subset V$ can be obtained by successive extensions of this kind.

 Notice that $\FF(V)$ can be identified with the subset of $\FF(U)\times\FF_x$ of all pairs that have the same restriction to $U\cap X_x$: the lemma is proven if we show that, for each $x\in X$, and $V\subset X_x$, the restriction $\FF(X_x)\frec\FF(V)$ is surjective. We proceed by induction on $\cork(x)$, where we define the {\em co-rank} of $x$, $\cork(x)$, as the maximum $r$ among the length of chains $x=x_0<\dots <x_r$, where $x_i\in X\;\forall i$.

If $\cork(x)=1$, the elements of $X_x$ have co-rank zero, so they are all open points of $X$: therefore $\FF(X_x)=\prod_{y\in X_x}\FF_y$, and for every $V\subset X_x$, $\FF(V)=\prod_{y\in V}\FF_y$. 

 Assume now that $\cork(x)>1$, and $V\subsetneqq X_x$. If $y\in X_x\smallsetminus V$,  $\FF(V\cup X_y)\cong \{(s,s')\in\FF(V)\times\FF_y\mid s_{|V\cap X_y}=s'_{|V\cap X_y}\}$. Given a section $s\in\FF(V)$, in order to extend it to $V\cup \bar{X}_y$ we need only to find some $s'\in \FF(X_y)$ such that $s'_{|V\cap X_y}=s_{|V\cap X_y}$, and we can do that by inductive hypothesis, as  $\cork(y)<\cork(x)$.  
\end{proof}


\subsubsection{Sheaves on fan spaces}
Let us now specialize to the case of fan spaces, that is, when $(X,\leq)$ is the set of cones of a regular fan $\Delta$, in $\R^n$. We introduce an order ``$\leq$'' on $\Delta$, defined in the following way:
$$
\sigma,\tau\in\Delta;\; \sigma\leq\tau \sse \tau\prec\sigma \quad\text{($\tau$ is a face of $\sigma$).}
$$  
Notice that $\leq$ {\em inverts} the face order, so the topology on $\Delta$ relative to this order depends only on the combinatorial structure of $\Delta$.

 Given a cone $\sigma\in\Delta$, the sets defined in the previous sections become: 
\begin{itemize}
\item[] $\bar{\Delta}_{\sigma}=\bar{\sigma}$ is the subfan generated by the cone $\sigma$;
\item[] $\Delta_{\sigma}=\bar{\Delta}_{\sigma}\smallsetminus \{\sigma\}=\bar{\sigma}\smallsetminus \sigma=\partial\sigma$ is the boundary of $\bar{\sigma}$;   
\item[] $\bar{\Delta}^{\sigma}=\{\tau\in\Delta:\sigma\prec\tau\}=\st \sigma$ is the \emph{star} of $\sigma$ in $\Delta$;
\item[] $\Delta^{\sigma}=\bar{\Delta}^{\sigma}\smallsetminus \{\sigma\}=\st\sigma\smallsetminus \sigma$.
\end{itemize}
  The open sets for this topology are the subfans of $\Delta$. The cone containing only the point $0\in \R^n$ is the unique maximal element, therefore the maximum, of $\Delta$. The set $\{0\}$ is open, and dense in $\Delta$, which is thus irreducible, as a topological space.  The maximal cones are the minimal elements with respect to the order. Notice that this topology corresponds to the Zariski topology on $X(\Delta)$ whose closed sets are the $T$-invariant subsets of $X(\Delta)$.

For an abelian group $G$, let $\tilde{G}$ be the constant presheaf with values in $G$. In fact $\tilde{G}$ is a sheaf, since all nonempty opens in $\Delta$ are connected. Moreover, $\tilde{G}$ is flabby, as all restrictions to nonempty opens are the identity of $G$. Also the sheaves, obtained by restricting $\tilde{G}$ to the stars of cones (which are closed in the fan space $\Delta$) are flabby. \vspace{5pt}

Let us define a class of sheaves, whose cohomology can be easily computed in terms of the simplicial cohomology of certain subsets of $\Delta$.
\begin{defin}
 If $G$ is an abelian group, and $\sigma\in\Delta$, the {\em simple} sheaf $G(\sigma)$ with support in $\sigma$ and values in $G$ is defined in the following way on the stalks: 
$$
G(\sigma)_{\tau}= \begin{cases} G &\textrm{if $\tau=\sigma$,}\\
0 & \textrm{if $\tau\not=\sigma$;}
\end{cases}
$$
and all restrictions are zero. 
\end{defin}
\begin{oss}\label{stelle}
We recall that the projection of  $\st\sigma$ on the quotient $N(\sigma)_{\R}$  is a fan, that we denote with $\overline{\st\sigma}$ (see Section \ref{orbite}). It can be identified with a subspace of $\Delta$: to every cone of $\overline{\st\sigma}$ we associate its preimage (through the projection) in $N_{\R}$.\vspace{5pt}
\end{oss}
The following two facts are crucial. (A less direct proof can be found in \cite{bac})
\begin{pro}\label{semplici}
(See \cite[Lemma 3.1]{bac}) Let $G$ be an abelian group,  and $\sigma\in\Delta$. 
\begin{itemize}
\item[(i)] The global sections of the simple sheaf $G(\sigma)$ are
\begin{displaymath}
\Gamma(\Delta,G(\sigma))=H^0(\Delta,G(\sigma))=\begin{cases}G\quad\text{if $\sigma$ is a maximal cone,}\\0\quad\text{otherwise.}\end{cases}
\end{displaymath}
\item[(ii)] If $i\geq 1$, then
 \begin{equation}
H^i(\Delta,G(\sigma))\cong \HR^{i-1}(S_{\overline{\st\sigma}},G),
\end{equation}
where $\HR^*$ is the reduced simplicial cohomology.
\end{itemize}
\end{pro}
\begin{proof}
The first part of the proposition is straightforward, in view of the formula \eqref{sezioni} describing the sections of sheaves on $\Delta$: the global sections of $G(\sigma)$ are the families of sections on maximal cones, compatible in the intersections. If $\sigma$ is not maximal, all these sections are zero. If it is maximal, the only nonzero sections are those in $G(\sigma)_{\sigma}=G$.

Let us prove (ii). As we pointed out in Remark \ref{stelle}, we can consider the star of  $\sigma$ as the closed image of a fan space $\overline{\st\sigma}$ through a continuous map $j_{\sigma}:\overline{\st\sigma}\iniett\Delta$ ($\sigma$ in $\Delta$  corrisponds to $0$ in $\overline{\st(\sigma)}$). The sheaf $G(\sigma)$ on $\Delta$ can be obtained by pushing forward via $j_{\sigma}$ the sheaf $G(0)$ on $\overline{\st\sigma}$: $G(\sigma)=(j_{\sigma})_* G(0)$. So $H^i(\Delta,G(\sigma))=H^i(\overline{\st\sigma},G(0))$. We must thus prove the follwing:
$$ \text{given a fan $\Delta$, }H^i(\Delta,G(0))=\HR^{i-1}(S_{\Delta},G) \text{ for $i>0$}.$$
Let $\tilde{G}$ the constant sheaf on $\Delta$, with values in $G$. $G(0)$ is a subsheaf of $\tilde{G}$. The quotient sheaf $\bar{G}=\tilde{G}/G(0)$ has stalks equal to those of $G$ everywhere, except at $0$, and, for any pair of cones $\tau\prec\sigma$, the restriction $\FF_{\sigma}\frec\FF_{\tau}$ is either zero or the identity. In other terms, if $i_0:\{0\}\iniett \Delta$, e $j_0:\Delta\smallsetminus 0$ are the (open and closed resp.) embeddings of $\{0\}$ and of its complement, then $\bar{G}=(j_0)_*(G_{|\Delta\smallsetminus 0})$. The short exact sequence
$$
0\frec G(0)=(i_0)_!\tilde{G}_{|0}\frec\tilde{G}\frec \bar{G}\frec 0
$$
induces the long exact sequence in cohomology:
\begin{multline}
0\frec H^0(\Delta,G(0))\frec H^0(\Delta,\tilde{G})\frec H^0(\Delta,\bar{G})\frec\\ \frec H^1(\Delta,G(0))\frec H^1(\Delta,\tilde{G})=0,
\end{multline}
and, since $\tilde{G}$ is acyclic,
$$
0\frec H^{i-1}(\Delta,\bar{G})\xrightarrow{\cong} H^i(\Delta,G(0))\frec 0\quad \forall\;i\geq 2.
$$
We can conclude by applying the following lemma.
\end{proof}
\begin{lemma}(See \cite[Theorem 2.1]{bac})\\
Given a fan $\Delta$ and an abelian group $G$, 
$$H^i(\Delta\smallsetminus 0,\bar{G})\cong\HR^i(S_{\Delta},G),$$ 
where $\bar{G}$ is the constant sheaf with values in $G$, on  $\Delta\smallsetminus 0$. 
\end{lemma}
\begin{proof}
Let  $\Delta_1=\{\rho_1,\dots,\rho_m\}$ be the set of 1-dimensional cones, and $\SSS$ the set of all the stars of the cones in $\Delta_1$: $\SSS=\{C_j=\st(\rho_j)\mid j=1,\dots,m\}$. $\SSS$ is a closed covering of $\Delta\smallsetminus 0$. The restriction of $\bar{G}$ to any intersection of elements in $\SSS$ is a flabby sheaf: we can thus consider a generalized Mayer-Vietoris sequence, that will allow us to compute the cohomology of $\bar{G}$. 

 Let $G_{i_0\dots i_p}=\bar{G}_{|C_{i_0}\cap\dots\cap C_{i_p}}$.  The sheaves appearing in the following complex are all acyclic, except (possibly) $\bar{G}$:
\begin{multline}\label{risoluz}
0\frec \bar{G}\frec \bigoplus_{i=1}^m G_i\frec \bigoplus_{1\leq i_0<i_1\leq m}G_{i_0i_1}\frec \bigoplus_{1\leq i_0<i_1<i_2\leq m}G_{i_0 i_1 i_2}\frec\dots\\
 \frec \bigoplus_{1\leq i_0<\dots<i_{n-1}\leq m}G_{i_0\dots i_{n-1}}\frec \bigoplus_{1\leq i_0<\dots<i_{n}\leq m}G_{i_0\dots i_{n}}\frec 0.  
\end{multline}
Differentials are defined as follows:
\begin{gather}d((a_{1_0\dots i_k})_{i_0\dots i_k})=(b_{j_0\dots j_{k+1}})_{j_0\dots j_{k+1}},
\end{gather}
where (indexes marked with $\wedge$ are omitted):
\begin{gather}\label{differ}
b_{j_0\dots j_{k+1}}=\sum_{h=0}^{k+1} (-1)^h \bar{a}_{j_0\dots\widehat{j_h}\dots j_{k+1}}.
\end{gather}
 The above notation means: if $a\in G_{j_0\dots \widehat{j_h}\dots j_{k+1}}(U)$, then $\bar{a}$ is the image of $a$ via the map (obvious definition) $G_{j_0\dots \widehat{j_h}\dots j_{k+1}}(U)\frec G_{j_0\dots  j_{k+1}}(U)$. 

The complex \eqref{risoluz} is exact (and so an acyclic resolution of the sheaf $\bar{G}$). To show this, fix a cone $\sigma$, say of dimension $d-1$, and consider the complex of the stalks relative to $\sigma$. Notice that $(\bar{G}_{i_0\dots i_p})_{\sigma}=G$ if $\sigma\in C_{i_0}\cap\dots\cap C_{i_p}$, and $0$ otherwise, so this localization is the Koszul complex $G\otimes \KKK^*(1,\dots,1)$:   
$$
0\frec G\frec G^d \frec \bigwedge^2 G^d\frec\dots\frec \bigwedge^d G^d\frec 0,
$$
which is exact.

 For any $(p+1)$-uple of indexes $i_0,\dots,i_p$, the intersection of the stars of 1-dimensional cones $\{ \rho_{i_k}\}$ is either the empty set, or the star of the cone $\sigma_{i_0\dots i_p}$ generated $\rho_{i_0},\dots \rho_{i_p}$:
$$
C_{i_0}\cap \dots\cap C_{i_p}=\begin{cases}\st{\sigma_{i_0\dots i_p}}& \text{if $\sigma_{i_0\dots i_p}\in\Delta$,}\\
\varnothing & \text{if $\sigma_{i_0\dots i_p}\not\in\Delta$.}\end{cases}
$$  
 Therefore, the complex of the global sections of resolution \eqref{risoluz} is 
\begin{multline}
0\frec \bigoplus_{\rho\in\Delta_1}H^0(\st(\rho),\bar{G}_{|\st(\sigma)})\frec \bigoplus_{\sigma\in\Delta_2}H^0(\st(\sigma),\bar{G}_{|\st(\sigma)})\frec\dots\\ \dots\frec \bigoplus_{\sigma\in\Delta_{n-1}}H^0(\st(\sigma),\bar{G}_{|\st(\sigma)})\frec\bigoplus_{\sigma\in\Delta_{n}}H^0(\{\sigma\},\bar{G}_{|\st(\sigma)})\frec 0.
\end{multline}
Since $H^0(\st(\sigma),\bar{G}_{|\st(\sigma)})=G$ for any $\sigma$, we can simplify as follows:
$$
0\frec \bigoplus_{\sigma\in\Delta_1}G_{(\sigma)}\frec\dots \frec\bigoplus_{\sigma\in\Delta_n}G_{(\sigma)}\frec 0.
$$ 
It is clear that, if we define  differentials as in \eqref{differ}, this is the cochain complex of $S_{\Delta}$ with values in $G$, and  its cohomology is the reduced simplicial cohomology of $S_{\Delta}$ with values in $G$.
\end{proof}
We define now a sheaf of $RT-$algebras on a fan space $\Delta$, which will play an important role in the following.
\begin{defin}
We denote with $\AAA$ the sheaf of $RT-$algebras on $\Delta$ with stalks:
$$
\AAA_{\sigma}=RT_{\sigma}\quad\sigma\in\Delta, 
$$ 
where $RT_{\sigma}\cong \Z\hat{T}$ is the ring of representations of the torus $T_{\sigma}$. The restrictions $\AAA_{\sigma}\frec \AAA_{\tau}$ are the maps $RT_{\sigma}\frec RT_{\tau}$ induced by the inclusions $\tau\iniett\sigma$.
\end{defin}
\begin{oss} Theorem \ref{immersione}, which gives an embedding of the equivariant K-theory of $X$ into a product of algebras of representations of tori, can be restated in a short (and useful) way, by means of the sheaf we have just defined:  
$$\Gamma(\Sigma,\AAA)=\AAA(\Sigma)\cong K_0^T(X(\Sigma))$$
for any subfan  $\Sigma\subseteq\Delta$.
\end{oss}
\begin{lemma}\label{Afiacco}
The sheaf $\AAA$ is flabby.
\end{lemma} 
\begin{proof}
(A similar construction is used in the proof of  \cite[3.2 Proposition]{bri})
Let us apply Lemma \ref{fiacc}: it is enough to show that, for any cone $\sigma\in\Delta$, the restriction to the boundary $\AAA_{\sigma}\frec\AAA(\partial\sigma)$ is surjective. Assume that $\sigma$ is maximal of dimension $n$ (we can always reduce to this case). As $\sigma$ is regular, it is generated by a basis  $\{u_1,\dots,u_n\}$ of the lattice $N\subset\R^n$. Let $\{v_1,\dots,v_n\}$ be a dual basis in $M$. For a face $\tau=\pos(u_{i_1},\dots,u_{i_d})\prec \sigma$, the quotient $M/\tau^{\perp}$ is generated (as a free abelian group) by the images of $v_{i_1},\dots,v_{i_d}$. So there is a section $\varphi_{\tau}: RT_{\tau}\cong\Z[M/\tau^{\perp}]\iniett RT \cong\Z[M]$, of the projection $\pi_{\tau}:\Z[M]\to \Z[M/\tau^{\perp}]$, that is, a map of $\Z$-algebras such that $\pi_{\tau}\varphi_{\tau}=id_{RT_{\tau}}$. Similarly, for every pairs of faces $\tau_1\prec\tau_2$ one finds injections $\varphi^{\tau_1}_{\tau_2}:RT_{\tau_1}\iniett RT_{\tau_2}$ such that $\pi^{\tau_2}_{\tau_1}\varphi^{\tau_1}_{\tau_2}=id_{RT_{\tau_1}}$, where $\pi^{\tau_2}_{\tau_1}$ is the projection $RT_{\tau_2}\to RT_{\tau_1}$. Now, an element of $\AAA(\partial\sigma)$ is represented by a family $(a_{\tau})_{\tau\prec\sigma}$, with $a_{\tau}\in RT_{\tau}$ such that ${a_{\tau_2}}_{|\tau_1}=a_{\tau_1}$, for $\tau_1\prec\tau_2\prec \sigma$.  Let $a\in RT$ be the element 
$$
a=\sum_{\tau\prec\sigma}(-1)^{n-\dim\tau+1}\varphi_{\tau}(a_{\tau}).   
$$
Then $a_{|\rho}=a_{\rho}$ for every $\rho\prec\sigma$. This follows easily by noticing that $\varphi_{\tau}(a_{\tau})_{|\rho}=\varphi^{\tau\cap\rho}_{\rho}({a_{\tau}}_{|\tau\cap\rho})$, and that 
$$\sum_{\rho\prec\tau}(-1)^{n-\dim\tau+1}=1,\quad\textrm{while, for $\rho'\precneqq\rho$,}\quad \sum_{\rho\cap\tau=\rho'}(-1)^{n-\dim\tau+1}=0.$$

\end{proof}


\subsection{$\Tor ^{RT}_p(K_0^T(X(\Delta)),\Z)$ and the reduced homology of $S_{\Delta}$}
\begin{teo}\label{principale}Let $X(\Delta)$ be the smooth toric variety associated to the fan $\Delta$ in $\R^n$. Assume that the maximal cones of $\Delta$ are all of maximal dimension: $\Delta_{\max}=\Delta_n$. Moreover assume that 
\begin{equation}\label{nullhom}
\HR_i(\lk\sigma,\Z)=0\qquad\forall\sigma\in S_{\Delta} \;,\forall i<\dim(\lk\sigma).
\end{equation}
 Then, for any $p>0$, 
\begin{equation}
\label{tor}
\Tor^{RT}_p(K_0^T(X(\Delta)),\Z)  \cong  \bigoplus_{i=p+1}^n\HR^{^{i-p-1}}\Bigl(S_{\Delta},\Z^{^{\binom{n}{i}}}\Bigr).
\end{equation}
\end{teo}
\begin{proof}
For convenience we split the proof in five steps.

\emph{Step 1.}
Let $\KKK^*$ be the complex of sheaves of $RT$-algebras on $\Delta$, obtained by tensoring with $-\otimes_{RT}\AAA$ the Koszul resolution of $\Z$ (as a $RT$-module: see \eqref{Koszul} in section \ref{RT}, setting $\sigma=0$). 

The terms of this complex are:
$$
\KKK^i={\bigwedge}^{-i} \AAA^n \quad\textrm{ for }-n\le i \le 0 
$$  
and zero for $i<-n$ and $i>0$.

Notice that, for $0\leq p\leq n$, the $(-p)$-th cohomology group of the complex $\KKK^*(\Delta)$ of the global sections of $\KKK^*$ is precisely $\Tor ^{RT}_{p}(K_0^T(X),\Z)$, since $\KKK^*(\Delta)$ is obtained by tensoring a projective resolution of $\Z$ with $K_0^T(X)$. On the other hand, the sheaves $\KKK^i$ are acyclic, as they are finite direct products of $\AAA$, which is flabby (Lemma \ref{Afiacco}). Hence for any $i\leq 0$, $H^i(\KKK^*(\Delta))$ is isomorphic to the $i$-th group of hypercohomolgy $\HHH^i(\KKK^*)$. So we have:
$$
\Tor ^{RT}_{p}(K_0^T(X),\Z)\cong\HHH^{-p}(\KKK^*)\quad\forall p\geq 0.
$$

\emph{Step 2.} 
  We will approximate $\HHH^*(\KKK^*)$ by means of a spectral sequence, relative to a suitable filtration of the complex $\KKK^*$ associated to the filtration of the fan space $\Delta$ by dimension of the cones. (This method is similar to that used in \cite[\S 4]{bac}, where a spectral sequence approximates the cohomology of a sheaf)

 Let $\Delta_i=\{\sigma\in\Delta\mid \dim(\sigma)=i\}$, $\codim(\sigma)=n-\dim(\sigma)$, and let $\{P^i\}_{i}$ be the filtration of $\Delta$, with $P^i$ the (closed) union of all cones of codimension $<i-1$:
\begin{multline*}
\varnothing=:P^1\subset \Delta_n=:P^2\subset\Delta_n\cup\Delta_{n-1}=:P^3\subset\dots\\\dots\subset\Delta_n\cup\dots\cup\Delta_1=\Delta\smallsetminus 0=:P^{n+1}\subset\Delta=:P^{n+2} 
\end{multline*}
Notice that 
\begin{align*}
P^{s+1}\smallsetminus P^s&=\{\sigma\in \Delta\mid\codim(\sigma)=s-1\}\:;\\
\Delta\smallsetminus P^s&=\{\sigma\in \Delta\mid\codim(\sigma)\geq s-1\}.
\end{align*}
Given a sheaf of abelian groups $\FF$ on $\Delta$, denote with $\FF_s$ the sheaf obtained by restricting $\FF$ to $\Delta\smallsetminus P^s$ and then by extending it by zero to $\Delta$. In other words, if $j_s:\Delta\smallsetminus P^s\iniett\Delta$, then $\FF_s=j_!(\FF_{|\Delta\smallsetminus P^s})$. The stalks of $\FF_s$ are:
$$
(\FF_s)_{\sigma}=\left\{ \begin{array}{ll} \FF_{\sigma}\quad &\textrm{if $\codim(\sigma)\geq s-1$}\\
0 & \textrm{otherwise.}
\end{array}\right.
$$    
For any $s$ there is an inclusion of sheaves $\FF_{s+1}\iniett\FF_s$. The support of the quotient sheaf $\FF_{s}/\FF_{s+1}$ is contained in $P^{s+1}\smallsetminus P^s$, which is a discrete topological space: so we have the following decomposition as a direct sum of simple sheaves:
\begin{equation}\label{somma}
\FF_{s}/\FF_{s+1}=\bigoplus_{\sigma\in\Delta_{n-s+1}}\FF_{\sigma}(\sigma).
\end{equation}
If, instead of a sheaf, we consider a {\em complex} of sheaves $\KKK$, we can define in a similar way the complexes $\KKK_s$, which give a filtration of $\KKK$, and the quotient complexes $\KKK_s/\KKK_{s+1}$, which can be written as a direct sums of complexes of simple sheaves.
 We can reconstruct the hypercohomology of the complex $\KKK$ from that of its restrictions to the sets $P^{s+1}\smallsetminus P^s$, by means of the following proposition, analogous to \cite[Thm 4.1]{bac}.
\begin{pro}
There is a cohomology spectral sequence:
\begin{equation}\label{hyperspettr}
E_1^{pq}=\HHH^{p+q}(\Delta,\KKK_{p}/\KKK_{p+1})\ergo \HHH^n(\Delta,\KKK).
\end{equation}
The differentials at level $E_1$ are the connecting homomorphisms of the long exact sequence relative to the short exact sequences
$$
0\fr\KKK_{p+1}/\KKK_{p+2}\fr\KKK_p/\KKK_{p+2}\fr\KKK_p/\KKK_{p+1}\fr 0.
$$
\end{pro} 
\begin{proof}
This is a quite standard result. Consider, for every $p$, the short exact sequence 
$$
0\to\KKK_{p+1}\to\KKK_p\to\KKK_p/\KKK_{p+1}\to 0.
$$
Let $\bar{\KKK}_p=\KKK_p/\KKK_{p+1}$. The hypercohomology long exact sequence is
$$
\dots\to\HHH^{i-1}(\bar{\KKK}_p)\stackrel{d}{\to}\HHH^i(\KKK_{p+1})\to\HHH^i(\KKK_{p})\to\HHH^{i}(\bar{\KKK}_p)\stackrel{d}{\to}\HHH^i(\KKK_{p+1})\to\dots
$$  
We set $A^{p,i}=\HHH^i(\KKK_{p})$ and $E^{p,i}=\HHH^{i}(\bar{\KKK}_p)$, change indexes, $q=i-p$, and then set $A=\oplus A^{pq}$, $E=\oplus E^{pq}$. We obtain an exact couple
$$
\xymatrix{A \ar[rr] &  & A \ar[dl]\\ 
            & E \ar[ul]^d &       }
$$
whose corresponding spectral sequence is the above one.
\end{proof}

\emph{Step 3.}
 Let us modify condition \eqref{nullhom} in order to simplify the spectral sequence just introduced. If $\lk\sigma$ is the link of $\sigma$ in the simplicial complex $S_{\Delta}$, we know that 
\begin{equation}
\HR_i(\lk\sigma,\Z)=0\qquad\forall\sigma\in\Delta\smallsetminus 0\;,\forall i<\dim(\lk\sigma).
\end{equation}
 All we can say about the highest degree homology of $\lk\sigma$ is that, if $d=\dim(\lk\sigma)$,
$$
\HR_d(\lk\sigma,\Z)\cong\left\{ \begin{array}{ll} \Z\quad &\textrm{if $\lk\sigma\cong S^d$ (i.e. $\sigma$ is ``inside'' $\Delta$)}\\
0 & \textrm{otherwise.}
\end{array}\right.
$$ 
This amounts to saying that the support of $S_{\Delta}$ is a (topological) manifold.

Applying the Universal Coefficient Theorem, we can replace homology with cohomology:
$$
\HR^i(\lk\sigma)\cong \HR_i(\lk\sigma)\quad\forall \sigma\;,\forall i
$$ 

By \eqref{somma} the term $E_1^{pq}$ of the spectral sequence \eqref{hyperspettr} can be written:
\begin{eqnarray*}
\HHH^{p+q}(\Delta,\KKK_p/\KKK_{p+1}) & = &\HHH^{p+q}\Bigl(\Delta,\bigoplus_{\sigma\in\Delta_{n-p+1}}\KKK_{\sigma}^*(\sigma)\Bigr)\\
 & = & \bigoplus_{\sigma\in\Delta_{n-p+1}}\HHH^{p+q}(\Delta,\KKK_{\sigma}^*(\sigma)).
\end{eqnarray*}
\emph{Step 4.} Every summand in the last direct sum is the hypercohomology group of a complex of simple sheaves on $\Delta$, with the same support $\sigma$. Let us study the behaviour of such complexes.

Let $G^*$ be a complex (in cohomology) of abelian groups, such that $G^i=0$ for $i\geq 0$. For any $\sigma\in\Delta$  there is a complex of simple sheaves $G^*(\sigma)$ supported in $\sigma$. Let us compute its hypercohomology. 
 
Consider the well-known spectral sequence:
$$
E_2^{ij}=H^i(\Delta,\HH^j(G^*(\sigma)))\ergo\HHH^{i+j}(\Delta,G^*(\sigma)),
$$
where $\HH^j$ denotes the $j$-th cohomology sheaf of a complex of sheaves. Here $\HH^j(G^*(\sigma))=(H^j(G^*))(\sigma)$, so by lemma \ref{semplici}, if $i>0$,
$$
H^i(\Delta,\HH^j(G^*(\sigma)))=\HR^{i-1}(\lk\sigma,H^j(G^*)).
$$
If $\sigma\not=0$, condition \eqref{nullhom} on $\Delta$ implies the vanishing of all columns $E_2^{i*}$, except possibly for $i=0$ and $i=n-\dim(\sigma)$. 

There are two cases:
\begin{enumerate}
\item $\dim(\sigma)=n$: only column $E_2^{0*}$ is nonzero, and $$\HHH^j(\Delta,G^*(\sigma))=E_{\infty}^{0j}=E_2^{0j}=H^0(\Delta,H^j(G^*)(\sigma))=H^j(G^*)$$ for $j\leq 0$;
\item $\dim(\sigma)<n$: then terms in the $0$-th column are
$$E_2^{0j}=H^0(\Delta,H^j(G^*)(\sigma))=0,$$
while those in the other non vanishing column, with index $n-\dim(\Delta)=d+1$, are
$$
E_2^{d+1,j}=H^{d+1}(\Delta,H^j(G^*)(\sigma)).
$$
\end{enumerate}
Now specialize to $G^*(\sigma)=\KKK_{\sigma}^*(\sigma)$: by formulae \eqref{TorRT} in \ref{RT}, the two cases can be rewritten:
\begin{enumerate}
\item $\dim(\sigma)=n$: $$\HHH^j(\Delta,\KKK^*(\sigma))=E_{\infty}^{0j}=E_2^{0j}=H^j(\KKK_{\sigma}^*)=\Tor ^{RT}_{-j}(RT_{\sigma},\Z)$$ 
so $\HHH^j(\Delta,\KKK^*(\sigma))=\Z$ for $j=0$, and $0$ for $j<0$;
\item $\dim(\sigma)<n$: in the $0$-th column we still have, for any $j$
$$E_2^{0j}=H^0(\Delta,H^j(\KKK_{\sigma}^*)(\sigma))=0,$$
while if the column index is $n-\dim(\sigma)=d+1$ we have
$$
E_2^{d+1,j}=H^{d+1}(\Delta,H^j(\KKK_{\sigma}^*)(\sigma))=H^{d+1}(\Delta,\Tor ^{RT}_j(RT_{\sigma},\Z)).
$$
The last term vanishes for  $j>n-\dim(\sigma)=d+1$: we conclude that 
$$
\HHH^j(\Delta,\KKK_{\sigma}^*(\sigma))=0
$$
for $j<0$. However, for $j\geq 0$ we can only say that 
$$
\HHH^j(\Delta,\KKK_{\sigma}(\sigma))=\HR^{d+1}(\lk\sigma,\Tor ^{RT}_{-j}(RT_{\sigma},\Z))=\HR^{d+1}(\lk\sigma,\wedge^{-j}\sigma^{\perp}).
$$
\end{enumerate}

\begin{center}
\begin{table}
\caption{Non-vanishing $E_{\infty}$ terms of the spectral sequence $E_1^{pq}=\HHH^{p+q}(\Delta,\KKK_p/\KKK_{p+1})\ergo \Tor ^{RT}_{-p-q}(K_0^T(X),\Z)$. (Empty cells contain zeroes, $*$ means ``unknown'')} 
\begin{tabular}{|l|c|c|c|c|c|c}\hline
   & $\scriptstyle{1}$   &  $\scriptstyle{2}$  & $\dots$ & $\scriptstyle{n}$   &$\scriptstyle{\mathbf{n+1}}$ &  $\scriptstyle{n+2}$ \\\hline
$\scriptstyle{-1}$ & $*$   &  $*$  & $\dots$ &  $*$  & $ *$  &  \\ \hline
$\scriptstyle{-2}$ &    &  $*$  &  $\dots$ &  $*$  &  $*$  & \\ \hline
$\scriptstyle{-3}$ &    &    &  $\dots$ &  $*$  &  $*$  &  \\ \hline
$\dots$ &$\dots$ &$\dots$ &$\dots$ &$\dots$ &$\dots$ &$\dots$ \\ \hline
$\scriptstyle{-n}$ &  &  & $\dots$ & $*$ & $*$ &  \\ \hline
$\scriptstyle{\mathbf{-n-1}}$&  &  & $\dots$ & & $*$ &  \\ \hline\rule{0pt}{12pt}
$\scriptstyle{-n-2}$&  &  & $\dots$ & & $\Tor^{RT}_1(K_0,\Z)$ & \\ \hline\rule{0pt}{12pt}
$\scriptstyle{-n-3}$&  &  & $\dots$ & & $\Tor^{RT}_2(K_0,\Z)$ & \\ \hline\rule{0pt}{12pt}
$\dots$&  &  & $\dots$ & & $\Tor^{RT}_3(K_0,\Z)$ & \\ \hline
$\dots$&$\dots$ &$\dots$ &$\dots$ &$\dots$ &$\dots$ &$\dots$ \\ \hline\rule{0pt}{12pt}
$\scriptstyle{-2n-1}$&  &  & \dots &  & $\Tor^{RT}_n(K_0,\Z)$ & \\ \hline
$\dots$&  &  & $\dots$ & &  &  \\
\end{tabular}
\end{table}
 \end{center}

\emph{Step 5.} In short: for $j<0$, the only contribution to $\HHH^j(\Delta,\KKK^*)$ comes from the $(n+1)$-th column  of the spectral sequence \eqref{hyperspettr},  relative to the restriction of the complex  $\KKK$ to the 0-dimensional cone. 

 So,  consider the case $\sigma=0$. The complex $\KKK_0^*(0)$, which is supported at the origin, has {\em  zero differentials}. In fact the differentials of $\KKK^*_0$ are $d\otimes id_{\Z}$, where $d$ are the differentials of a Koszul complex relative to elements that belong to the annihilator of $\Z$ as $RT$-module. 

Therefore $\KKK_0^*(0)$ is a direct sum of the complexes $\KKK^*_{0i}$, each one having at most a nonzero component:
$$
\KKK_0^*(0)=\bigoplus_{i=-n}^0 \KKK^*_{0i}(0)
$$
where we define, for $-n\leq m\leq 0$,
$$
\KKK^m_{0i}(0):=\begin{cases} \KKK^i_0(0) & \textrm{if $m=i$,}\\ 0 & \textrm{if $m\not=i$.}
 \end{cases}
$$ 
Hypercohomology commutes with direct sums, and the hypercohomology of a complex with a unique non-vanishing sheaf $\FF$ in degree $i$ coincides (up to a shifting of the indexes by $i$) with the cohomology of $\FF$. So:   
$$
\HHH^l(\Delta,\KKK^*_{0i}(0))=H^{l-i}(\Delta,\KKK^i_{0i}(0))\quad\text{for $l\geq i$}.
$$
Since $\KKK^{-i}_0(0)=\Tor^{RT}_i(\Z,\Z)=\Z^{^{\binom{n}{i}}}$, we conclude that, {\em for $p>0$}, 
\begin{eqnarray}
\notag
\Tor ^{RT}_p(K_0^T(X),\Z) & = & \bigoplus_{i=p}^{n}H^{i-p}(\Delta,\KKK^{-i}_0(0))\\ & = & \bigoplus_{i=p+1}^n\HR^{i-p-1}\Bigl(S_{\Delta},\Z^{^{\binom{n}{i}}}\Bigr).
\end{eqnarray}
And the theorem is proved.  
\end{proof}
\begin{oss}
What can be said about $\Tor ^{RT}_0(K_0^T(X),\Z)\cong\Z\otimes_{RT} K_0^T(X)$? Since some terms $E_1^{ij}$ do not vanish for $i+j\geq 0$, differentials in $E_1$ or in higher levels forbids us to describe $\Z\otimes_{RT} K_0^T(X)$ as we did for higher $\Tor^{RT}$. However we can give an upper bound of its rank: it is less or equal to that of $\oplus_{p+q=0} E_1^{pq}$, since, for $p+q<0$ and $p<n+1$, we have $E_1^{pq}=0$, so $E_{\infty}^{pp}\subset E_1^{pp}$ for any $p$.  
\end{oss}
\subsubsection{Computing $\Tor_p^{RT}(K^T_q(X),\Z)$ for $q>0$}
We can find an expression of all the other columns of Merkurjev's spectral sequence, starting from \eqref{tor}. 

As we have already pointed out, for any $p>0$ there is a splitting short exact sequence (see \eqref{kappasup} in the proof of Prop.~\ref{piattezza}):
\begin{multline}
0\frec \Tor^{RT}_p(K_0^T(X),\Z)\otimes_{\Z} K_q(k)\frec \Tor^{RT}_p(K^T_q(X),\Z)\frec \\ \frec \Tor_1^{\Z}(\Tor_{p-1}^{RT}(K_0^T(X),\Z),K_q(k))\frec 0.
\end{multline}
Equivalently:
\begin{multline}\label{pappa}
\Tor^{RT}_p(K^T_q(X),\Z)\cong\Tor^{RT}_p(K_0^T(X),\Z)\otimes_{\Z} K_q(k)\oplus\\\Tor_1^{\Z}(\Tor_{p-1}^{RT}(K_0^T(X),\Z),K_q(k)).
\end{multline}
Notice that for $p=1$, since $\Z\otimes_{RT}K_0^T(X)$ is torsion free, the second summand on the right vanishes. Therefore all the terms $\Tor_*^{RT}(K_0^T(X),\Z)$ appearing in \eqref{pappa} can be replaced by the expression we gave in \eqref{tor}: we obtain expressions for any $E_{pq}^2$ term of Merkurjev's spectral sequence, involving only the reduced cohomology of $S_{\Delta}$ and the K-theory of the field $k$.

\subsection{An alternative proof of Proposition \ref{piattezza}}\label{alt}
Theorem \ref{principale} allows us to prove Proposition \ref{piattezza} without making reference to Reisner's Theorem (Proposition \ref{reisn}). One implication, $(4)\ergo(1)$ (here we assume $(1)\sse(2)\sse(3)$), follows immediately from Theorem \ref{principale}. The opposite one is not straightforward, but it can be easily seen after a couple of lemmas on Stanley-Reisner rings.
\begin{lemma}\label{local}Let $\Sigma$ be a simplicial complex on the vertices $\{v_1,\dots,v_t\}$, $\sigma=[v_1,\dots,v_l]$ a face of $\sigma$.   
Then we have an isomorphism of localizations of Stanley-Reisner rings:
$$
\Z[\Sigma]_{(X_{\sigma})}\cong\Z[\st\sigma]_{(X_{\sigma})},
$$ 
where $X_{\sigma}$ is the image of the monomial $X_1\dots X_l$.
\end{lemma}
\begin{proof}
Let the vertices of the star $\st\sigma$ be $\{v_1,\dots,v_r\}$ ($l\le r\le t$). By definition $\Z[\Sigma]=\Z[X_1,\dots,X_t]/I_{\Sigma}$, where $I_{\Sigma}=(X_{i_1}\dots X_{i_k}\mid [v_{i_1},\dots,v_{i_k}]\not\in\Sigma)$, while $\Z[\st\sigma]=\Z[X_1,\dots,X_r]/I_{\sigma}$, where $I_{\sigma}=(X_{i_1}\dots X_{i_k}\mid [v_{i_1},\dots,v_{i_k}]\not\in\st\sigma)$. After we localize to the multiplicative system $\{X_{\sigma}^n\mid n\geq 0\}$, all monomials $X_{i_1}\dots X_{i_k}\in\Z[\Sigma]$ such that $[v_{i_1},\dots,v_{i_k}]\not\in \st\sigma$ vanish: so the inclusion $\Z[X_1,\dots,X_r]\iniett\Z[X_1,\dots,X_t]$ can be lifted to a well defined ring homomorphism $\Z[\Sigma]_{(X_{\sigma})}\to\Z[\st\sigma]_{(X_{\sigma})}$ that is easily seen to be injective and surjective.  
\end{proof}   
\begin{lemma}\label{starlink}
With the same notation as above, 
$$
\Z[\lk\sigma]=\frac{\Z[\st\sigma]}{(X_{\sigma})}.
$$ 
\end{lemma}
\begin{proof}
$
\Z[\lk\sigma]=\Z[X_{l+1},\dots,X_r]/I_{\lk\sigma}$, with $I_{\lk\sigma}=I_{\sigma}\cap\Z[X_{l+1},\dots,X_r]$.
\end{proof}
The two following propositions can replace the reference to Reisner's Theorem in the proof of Proposition \ref{piattezza}.
 
Remember that the closure of the orbit $\OO_{\sigma}$ can be considered a toric variety for the action of the quotient torus $\bar{T}_{\sigma}=T/T_{\sigma}$, and the simplicial complex associated to its fan is $\lk\sigma$ (see section \ref{orbite}).
\begin{pro}
Let $X=X(\Delta)$ be a smooth toric variety associated to the fan $\Delta$, pure of dimension  $\dim(X)$, and $Y=\overline{\OO}_{\sigma}$ the closure of the orbit associated to the cone $\sigma$. 
 Then
$$
K_0^{\bar{T}_{\sigma}}(Y)\textrm{ is not $R\bar{T}_{\sigma}$-flat }\quad\ergo\quad K_0^{T}(X)\textrm{ is not $RT$-flat}. 
$$   
\end{pro}
\begin{proof}
Since $K_0^{\bar{T}_{\sigma}}(Y)$ is not $R\bar{T}_{\sigma}$-flat, by lemma \ref{ohoh} (in the proof of $(4)\ergo(1)$ of Proposition \ref{piattezza}) $\Z[\lk\sigma]$ is not a Cohen-Macaulay ring. By Lemma \ref{starlink}, $\Z[\st\sigma]$ is not Cohen-Macaulay either, for $X_{\sigma}$ is a regular element. The localization map $\Z[\st\sigma]\to\Z[\st\sigma]_{(X_{\sigma})}$ is flat, so $\Z[\st\sigma]_{(X_{\sigma})}$ can not be Cohen-Macaulay (\cite[Theorem 2.1.7]{bh}). By Lemma \ref{local} we conclude that $\Z[S_{\Delta}]$ is not Cohen-Macaulay, so $K_0^{T}(X)$ is not $RT$-flat.
\end{proof}
\begin{pro}
Let $X=X(\Delta)$ be a smooth toric variety associated to the fan $\Delta$, pure of dimension  $\dim(X)$, and suppose that, for some $\tau\in S_{\Delta}$ and some $i<\dim(\lk\tau)$, $\HR_i(\lk\tau)\not=0$. Then there is some orbit closure $Y=\overline{\OO}_{\sigma}$ such that $K_0^{\bar{T}_{\sigma}}(Y)$ is not $R\bar{T}_{\sigma}$-flat.  
\end{pro}
\begin{proof}
If (a) $Y$ satisfies the hypothesis \eqref{nullhom} of Theorem \ref{principale}, that is, if for all $\rho\in\lk\tau$ we have $\HR_i(\lk_{\rho}(\lk\tau))=0 \quad\forall\,i<\dim\lk_{\rho}(\lk\tau)$, then we can apply Theorem \ref{principale} to $\overline{\OO}_{\tau}$, viewed as a toric variety whose simplicial complex is $\lk\tau$, and conclude with $\sigma=\tau$. In fact the formula \eqref{tor} implies that $\Tor^{R\bar{T}_{\tau}}_1(K_0^{\bar{T}_{\tau}}(Y))\not=0$. 

If, on the contrary, (b) there is a $\rho$ such that $\HR_i(\lk_{\rho}(\lk\tau))\not=0 \quad\exists\,i<\dim\lk_{\rho}(\lk\tau)$, then we cannot apply Theorem \ref{principale}, but we can reason by induction, replacing $X$ with $\overline{\OO}_{\tau}$ and $\tau$ with $\rho$. Since $\dim\overline{\OO}_{\tau}<\dim X$, after a finite number of steps assumption (a) must be satisfied. We can conclude, by noting that if $\tau\prec\sigma\in\Delta$, $\OO_{\tau}$ is also a $\bar{T}_{\sigma}$ orbit, and its closure in $\overline{\OO}_{\sigma}$ is the same as its closure in $X$.
\end{proof}

\end{document}